\documentclass[twocolumn]{autart}    
\usepackage[square,comma,sort&compress]{natbib}
\usepackage[colorlinks=true,bookmarksopen,bookmarksnumbered,citecolor=blue,urlcolor=blue]{hyperref}

\usepackage{mathrsfs}
\usepackage{graphics}
\usepackage{amssymb}
\usepackage{amsmath}
\usepackage{color}
\usepackage{xspace}
\usepackage{algpseudocode}
\usepackage{bbm}
\usepackage{comment}
\usepackage{algorithmicx}
\usepackage{subfig}
\usepackage{psfrag}
\usepackage{soul}
\usepackage{tikz}
\usetikzlibrary{shapes,arrows}
\usetikzlibrary{positioning}

\tikzstyle{block}=[draw opacity=0.7,line width=1.4cm]

\definecolor{CranJ}{cmyk}{0,0.69,0.54,0.04} %cranberry jello
\definecolor{PinkJ}{cmyk}{0,0.71,0.43,0.12} %pink jeep
\definecolor{Cran}{cmyk}{0,0.73,0.41,0.29} %cranberry 
\definecolor{VRed}{cmyk}{0,0.75,0.25,0.2} %violetred
\definecolor{ORed}{cmyk}{0,0.75,0.75,0} %orangered4
\definecolor{CBlue}{cmyk}{1,0.25,0,0} %curacao	
\newcommand{\VV}{\mathcal{V}}
\newcommand{\EE}{\mathcal{E}}
\newcommand{\GG}{\mathcal{G}}

\newcommand{\lL}{\vect{L}}
\newcommand{\UHlambda}{\underline{\hat{\lambda}}}
\newcommand{\Ulambda}{\underline{\lambda}}
\newcommand{\Bdout}{\bar{\mathsf{d}}_{\operatorname{out}}}
\newcommand{\BnL}{\overline{\|\lL\|}}
\newcommand{\R}{\vect{\mathfrak{R}}}
\newcommand{\rr}{\vect{\mathfrak{r}}}

\newcommand{\PPi}{\vect{\Pi}}
\newcommand{\pPi}{\vect{\Pi}}
\newcommand{\real}{{\mathbb{R}}} \newcommand{\reals}{{\mathbb{R}}}
\newcommand{\realpositive}{{\mathbb{R}}_{>0}}
\newcommand{\realnonnegative}{{\mathbb{R}}_{\ge 0}}
\newcommand{\integersnonnegative}{{\mathbb{Z}}_{\ge 0}}

 \newcommand{\eps}{\epsilon}

\newcommand{\Hlambda}{\hat{\lambda}}

\newcommand{\map}[3]{#1:#2 \rightarrow #3}

\newcommand{\din}{\mathsf{d}_{\operatorname{in}}}
\newcommand{\dout}{\mathsf{d}_{\operatorname{out}}}

\newcommand{\until}[1]{\in\{1,\dots,#1\}}

\newcommand{\setdef}[2]{\{#1 \; |\; #2\}}

\newcommand{\argmax}{\operatorname{argmax}}
\newcommand{\vect}[1]{\boldsymbol{\mathbf{#1}}}
\newcommand{\Bvect}[1]{\bar{\boldsymbol{\mathbf{#1}}}}
\newcommand{\Tvect}[1]{\tilde{\boldsymbol{\mathbf{#1}}}}
\newcommand{\Hvect}[1]{\hat{\boldsymbol{\mathbf{#1}}}}

\newcommand{\dvect}[1]{\dot{\vect{#1}}}

\newcommand{\Sym}[1]{\operatorname{Sym}(#1)}
\newcommand{\Diag}[1]{\operatorname{Diag}(#1)}

\newcommand{\avrg}[1]{\frac{1}{N}\sum\nolimits_{j=1}^N\!\!#1^j}
\newcommand{\SUM}[2]{\sum\nolimits_{#1}^{#2}}
 \newcommand{\boxend}{\hfill \ensuremath{\Box}}

\newtheorem{theorem}{Theorem}[section]

\newtheorem{remark}[theorem]{Remark}

\newtheorem{proposition}[theorem]{Proposition}  

\newcommand{\oprocendsymbol}{\hbox{$\bullet$}}
\newcommand{\oprocend}{\relax\ifmmode\else\unskip\hfill\fi\oprocendsymbol}

\makeatletter
\renewcommand*{\@opargbegintheorem}[3]{\trivlist
      \item[\hskip \labelsep{ #1\ #2}] (#3):\ \itshape}
\makeatother

\begin{document}

\begin{frontmatter}

  \title{Distributed event-triggered communication for dynamic average
    consensus in networked systems\thanksref{footnoteinfo}}

  % Conference paper title: "Dynamic Average Consensus
  % with Distributed Event-triggered Communication"

  \thanks[footnoteinfo]{A preliminary version appears at the IEEE
    Conference on Decision and Control
    as~\citep{SSK-JC-SM:14-cdc}.} % Corresponding author: S. S. Kia}
  
  \author[Paestum]{Solmaz S. Kia}\ead{solmaz@uci.edu} \quad 
    \author[Rome]{Jorge Cort\'es}\ead{cortes@ucsd.edu}
  \quad 
  \author[Rome]{Sonia Mart{\'\i}nez}\ead{soniamd@ucsd.edu}
  
 \address[Paestum]{Department of Mechanical and Aerospace Engineering,
    University of California at Irvine, Irvine, CA 92697, USA}
  \address[Rome]{Department of Mechanical and Aerospace Engineering,
    University of California at San Diego, La Jolla, CA 92093, USA}
  
  \begin{abstract}
    This paper presents distributed algorithmic solutions that employ
    opportunistic inter-agent communication to achieve dynamic average
    consensus.   In our solutions each agent is endowed with a local criterion that enables it to determine whether to broadcast its state to  its neighbors.
    Our starting point is a continuous-time distributed coordination strategy that,
    under continuous-time communication, achieves practical asymptotic
    tracking of the dynamic average of the time-varying agents'
    reference inputs.  Then, for this algorithm, depending on the directed or undirected nature
    of the time-varying interactions and under suitable connectivity
    conditions, we propose two different distributed event-triggered
    communication laws that prescribe agent communications at discrete time
    instants in an opportunistic fashion.  In both cases, we establish
    positive lower bounds on the inter-event times of each agent and
    characterize their dependence on the algorithm design
    parameters. This analysis allows us to rule out the presence of
    Zeno behavior and characterize the asymptotic correctness of the
    resulting implementations.  Several simulations illustrate the
    results.
  \end{abstract}
\end{frontmatter}

\section{Introduction}

The dynamic average consensus problem consists of designing a
distributed algorithm that allows a group of agents to track the
average of individual time-varying reference inputs, one per agent.  This
problem has applications in numerous areas that involve distributed
sensing and filtering, including distributed
tracking~\citep{PY-RAF-KML:07}, multi-robot
coordination~\citep{PY-RAF-KML:08}, sensor
fusion~\citep{ROS-JSS:05,ROS:07}, and distributed
estimation~\citep{AC-TO-MT-RC-LS:13}. Our goal here is to develop
algorithmic solutions to the dynamic average consensus problem which
rely on agents %autonomously 
locally deciding when to share information with
their neighbors in an \emph{opportunistic} fashion for greater efficiency and
energy savings. By opportunistic, we mean that
the transmission of information to the neighbors should happen at
times when it is needed to preserve the stability and convergence of
the coordination algorithm. 

\emph{Literature review}: The literature of cooperative control has
proposed dynamic average consensus algorithms that are executed either
in
continuous-time~\citep{DPS-ROS-RMM:05b,ROS-JSS:05,RAF-PY-KML:06,HB-RAF-KML:10,
  SSK-JC-SM:14-ijrnc} or in fixed stepsize
discrete-time~\citep{MZ-SM:08a,SSK-JC-SM:14-ijrnc}.  Continuous-time
algorithms operate under the assumption of continuous agent-to-agent
information sharing. Although discrete-time algorithms are more
amenable to practical implementation, use of a fixed communication
step-size, which should be designed to address also rarely occurring
worst-case situations, can be a wasteful use of the network
resources. In addition, in these discrete-time algorithms
communication and computation stepsizes are tied together, resulting
in potentially a conservative stepsize for communication times. This
can result in costly implementations, as performing communication
usually requires more energy than computation. Moreover, the
assumption of periodic, synchronous communication is unrealistic in
many scenarios involving cyber-physical systems, as processors are
subject to natural delays and errors which may deviate them from the
perfect operational conditions the strategies are designed for.
Event-triggered communication offers a way to address these
shortcomings by prescribing in an opportunistic way the times for
information sharing and allowing individual agents to determine these
autonomously.  In recent years, an increasing body of work that seeks
to trade computation and decision-making for less communication,
sensing or actuation effort while guaranteeing a desired level of
performance has emerged, see
e.g.,~\citep{WPMHH-KHJ-PT:12,MMJ-PT:11,XW-MDL:11}. Closest to the
problem considered here are the works that study event-triggered
communication laws for static average consensus, see
e.g.,~\citep{DVD-EF-KHJ:12,EG-YC-HY-PA-DC:13,GSS-DVD-KHJ:13,YF-GF-YW-CS:13,CN-JC:14-acc}
and references therein.

\emph{Statement of contributions}: We propose novel algorithmic
solutions to the dynamic average consensus problem that employ
opportunistic strategies to determine the communication times among
neighboring agents.  The basic idea is that agents share their
information with neighbors when the uncertainty in the outdated
information is such that the monotonic convergent behavior of the
overall network can no longer be guaranteed. To realize this concept,
depending on the connectivity properties of the interaction topology,
we propose and characterize the asymptotic correctness of two
different distributed event-triggered communication laws. Our least
stringent connectivity conditions are modeled by a time-varying,
weight-balanced piecewise constant digraph which is jointly strongly
connected over an infinite sequence of contiguous and uniformly
bounded time intervals. In the second scenario, we consider
interaction topologies modeled by a time-varying, piecewise continuous
undirected connected graph, which allows us to further refine our
analytical guarantees.  By establishing positive lower bounds on the
inter-event times of each agent for both cases, we also show that the
proposed distributed event-triggered communication laws are free from
Zeno behavior (the undesirable situation where an infinite number of
communication rounds are triggered in a finite amount of
time). Finally, we analyze the dependence of the inter-event times on
the algorithm design parameters.  This characterization provides
guidelines on the trade-offs between the minimum inter-event times for
communication and the performance and energy efficiency of the
proposed algorithms.  We demonstrate through several comparative
simulation studies the advantages of our proposed event-triggered
communication strategies over schemes that rely on continuous-time
communication as well as discrete-time communication with fixed
stepsize.
    
% \emph{Organization}: Section~\ref{sec::prelim} gathers basic
% notation and concepts from graph
% theory. Section~\ref{se:problem-statement} presents the network
% model and the dynamic average consensus
% problem. Section~\ref{sec::Discrete} introduces our continuous-time
% algorithmic solutions with event-triggered communication.
% Section~\ref{sec::num} presents simulation results and
% Section~\ref{sec::conclu} gathers our conclusions and ideas for
% future work.

\section{Notation and terminology}
% In this section, we introduce basic notation and concepts from graph
% theory.
We let $\reals$, $\realnonnegative$, $\realpositive$, and
$\integersnonnegative$ denote the set of real, nonnegative real,
positive real, and nonnegative integer, respectively.  The transpose
of a matrix $\vect{A}$ is~$\vect{A}^\top$.  We let $\vect{1}_n$
(resp. $\vect{0}_{n}$) denote the vector of $n$ ones (resp. $n$
zeros). We let $\PPi_n = \vect{I}_n - \frac{1}{n}
\vect{1}_n\vect{1}_n^\top$, where $\vect{I}_n$ is the $n\times n$
identity matrix.  When clear from the context, we do not specify the
matrix dimensions.  For $\vect{u}\in\reals^d$,
$\|\vect{u}\|=\sqrt{\vect{u}^\top\vect{u}}$ is the standard Euclidean
norm. For $u \in \reals$, $|u|$ is its absolute value. For a
time-varying measurable locally essentially bounded signal $\vect{u} :
\realnonnegative \to\real^{m}$, we denote by
$\|\vect{u}\|_{\text{ess}}$
% = (\text{ess})\sup \setdef{\|\vect{u}(t)\|}{t \in \realnonnegative}$
the essential supremum norm. For a scalar signal $u$, we use
$|u|_{\text{ess}}$ instead.  For vectors
$\vect{u}_1,\dots,\vect{u}_m$, we let $\vect{u} =
(\vect{u}_1,\dots,\vect{u}_m)$ represent the aggregated vector. 
%For finite sets $V_1$ and $V_2$, we denote by $V_1\backslash V_2$ the set
%whose elements consist of all the elements of $V_1$ that are not in
%$V_2$. 
In~a networked system, we distinguish the local variables at
each agent by a superscript, e.g., $\vect{x}^i$ is the local state~of
agent $i$. If $\vect{p}^i\in\reals^d$ is a variable of agent $i$, the
aggregate of a network with $N$ agents is $\vect{p} =
(\vect{p}^1,\dots,\vect{p}^N) \in (\reals^d)^N$.  %To deal with the
%dynamics that follow, it is convenient to introduce the orthonormal  matrices
%$\vect{\mathfrak{T}} \in \real^{N \times N}$ % and ${\R} \in \real^{N \times
 % N-1}$, and the vector $\vect{r} \in \real^N$ as follows,
%Notice that $\vect{\mathfrak{T}}^\top\vect{\mathfrak{T}}=\vect{I}_N$.

%\subsection{Graph theory}\label{sec::graph}

\emph{Graph theory}: Here, we briefly review some basic concepts from
graph theory and linear algebra following~\citep{FB-JC-SM:09}.  A
\emph{directed graph}, or simply a \emph{digraph}, is a pair $\GG =
(\VV ,\EE )$, where $\VV=\{1,\dots,N\}$ is the \emph{node set} and
$\EE \subseteq \VV\times \VV$ is the \emph{edge set}.  For an edge
$(i,j) \in\EE$, $i$ is called an \emph{in-neighbor} of $j$ and $j$ is
called an \emph{out-neighbor} of~$i$.  We let $\mathcal{N}^i$ denote
the set of out-neighbors of~$i\in{\VV}$. A graph is \emph{undirected}
if $(i,j) \in \EE$ when $(j,i)\in\EE$.  A \emph{directed path} is a
sequence of nodes connected by edges.
% A \emph{directed path} is an ordered sequence of vertices such that
% any ordered pair of vertices appearing consecutively is an edge of
% the digraph.
A digraph is called \emph{strongly connected} if for every pair of
vertices there is a directed path connecting~them. Given digraphs
$\GG_i=(\VV,\EE_i)$, $i \until{m}$, their \emph{union} is the
graph $\cup_{i=1}^n \GG_i=(\VV,\EE_1\cup\EE_2\cup \dots\cup\EE_m)$.

A \emph{weighted digraph} is a triplet $(\VV ,\EE, \vect{\sf{A}})$,
where $\GG=(\VV ,\EE )$ is a digraph and $
\vect{\sf{A}}\in\real^{N\times N}$ is a weighted \emph{adjacency}
matrix with the property that $ \mathsf{a}_{ij} >0$ if $(i, j) \in\EE$
and $\mathsf{a}_{ij} = 0$, otherwise.
% We use $\Gamma(\vect{\sf{A}})$ to denote a digraph induced by a
% given adjacency matrix $\vect{\sf{A}}$.
A weighted digraph is \emph{undirected} if $\mathsf{a}_{ij} =
\mathsf{a}_{ji}$ for all $i,j\in\VV$. We refer to a strongly connected
and undirected graph as \emph{connected}. The \emph{weighted out-} and
\emph{in-degrees} of a node $i$ are, respectively, $\dout^i =\sum^N_{j
  =1} \mathsf{a}_{ij}$ and $\din^i =\sum^N_{j =1}
\mathsf{a}_{ji}$. %We let
% $\mathsf{d}^{\max}_{\text{out}}
% =\max\{\mathsf{d}_{\text{out}}^1,\cdots,\mathsf{d}_{\text{out}}^N\}$
% and $\mathsf{d}_{\min}^{\text{out}} = \underset{i \in
%   \until{N}}{\min} \mathsf{d}^{\text{out}} (i)$ denote the maximum
% and minimum weighted out-degree.
A digraph is \emph{weight-balanced} if at each node $i\in\VV$, the
weighted out- and in-degrees coincide (although they might be
different across different nodes).  The \emph{(out-) Laplacian} matrix
is $\lL = \vect{\mathsf{D}}^{\text{out}} - \vect{\mathsf{A}}$, where
$\vect{\mathsf{D}}^{\text{out}} = \Diag{\dout^1,\cdots, \dout^N} \in
\reals^{N \times N}$.  Note that $\lL\vect{1}_N=\vect{0}$. A digraph
is weight-balanced if and only if $\vect{1}_N^T\lL=\vect{0}$ if and
only if $\Sym{\lL} = \tfrac{1}{2} (\lL + \lL^T)$ is positive
semi-definite. Based on the structure of $\lL$, at least one of the
eigenvalues of $\lL$, denoted by $\lambda_1, \dots,\lambda_N$, is zero
and the rest of them have nonnegative real parts. We let $\lambda_1=0$
and $\Re(\lambda_i)\leq \Re(\lambda_j)$, for $i<j$, where $\Re(\cdot)$
denotes the real part of a complex number.
%We denote the eigenvalues of $\lL$ by $\lambda_1, \dots,\lambda_N$,
%where $\lambda_1=0$ and $\Re(\lambda_i)\leq \Re(\lambda_j)$, for
%$i<j$, and 
We denote the eigenvalues of $\Sym{\lL} $ by
$\hat{\lambda}_1,\dots,\hat{\lambda}_N$. For a strongly connected and
weight-balanced digraph, zero is a simple eigenvalue of both
$\vect{L}$ and $\Sym{\lL} $. In this case, we order the eigenvalues of
$\Sym{\lL} $ as $0 = \hat{\lambda}_1
<\hat{\lambda}_2\leq\hat{\lambda}_3 \leq \hdots\leq\hat{\lambda}_N$.

Throughout the paper, we deal with time-varying digraphs with fixed
vertex set.  A time-varying digraph $\realnonnegative\ni t \mapsto
\GG(t) = (\VV,\EE(t),\vect{\mathsf{A}}(t))$ is \emph{piecewise
  continuous (respectively, piecewise constant)} if the map $t\mapsto
\vect{\mathsf{A}}(t)$ is piecewise continuous (respectively, piecewise
constant) from the right.  In such case, we denote by $\{s_k\}_{k \in
  \integersnonnegative}$ the time instants at which this map is
discontinuous and refer to them as \emph{switching times}. By
convention, $s_0=0$.  A time-varying digraph $\realnonnegative\ni t
\mapsto \GG(t)$ has \emph{uniformly bounded weights} if, for all $t
\in \realnonnegative$, $0<\underline{\mathsf{a}}
\leq\mathsf{a}_{ij}(t) \in[\underline{\mathsf{a}}, \bar{\mathsf{a}}]$,
with $0<\underline{\mathsf{a}} \leq \bar{\mathsf{a}}$, if
$(j,i)\in\EE(t)$, and $\mathsf{a}_{ij}=0$ otherwise.  A time-varying
digraph $\realnonnegative\ni t \mapsto \GG(t)$ is \emph{strongly
  connected} if each $\GG(t)$ is strongly connected, and is
\emph{jointly strongly connected} over $[t_1,t_2)$ if ${\cup}_{t \in
  [t_1,t_2)}\GG(t)$ is strongly connected.  A piecewise constant
time-varying digraph $\realnonnegative\ni t \mapsto \GG(t)$ is
\emph{recurrently jointly strongly connected} if the sequence of
inter-switching times $\{s_{k+1}-s_k\}_{k \in \integersnonnegative}$
is uniformly lower bounded and there exists an infinite sequence of
contiguous uniformly bounded intervals $\{[s_{k_j},s_{k_{j+1}})\}_{j
  \in \integersnonnegative}$, with $s_{k_0}=s_0$, such that $\GG(t)$
is jointly strongly connected over $[s_{k_j},s_{k_{j+1}})$, for all $j
\in \integersnonnegative$.  Finally, a time-varying digraph
$\realnonnegative\ni t \mapsto \GG(t)$ is \emph{weight-balanced} if
each $\GG(t)$ is weight-balanced.  
% Throughout the paper, $\{\GG\}_{\{s_k\}}$ denotes a piecewise
% constant recurrently jointly strongly connected and weight-balanced
% digraph with uniformly bounded weights.
For a piecewise constant recurrently jointly strongly connected and
weight-balanced digraph with uniformly bounded weights, we let
$\{s_k^i\}_{k \in \integersnonnegative}\!\!\subseteq\!\{s_k\}_{k \in
  \integersnonnegative}$ be the times when an agent $i\!\in\!\!\VV$
acquires a new in-neighbor and we define
\begin{align*}
  ~\BnL & \!=\! \sup \setdef{\|\lL_t\|}{t\! \in\! \realnonnegative},
  ~~\UHlambda_2 \!=\! \inf \setdef{\Hlambda_2(\lL_t)}{t\! \in\!
    \realnonnegative} ,
  \\
  ~{\Bdout}^{i} & \!=\! \sup \setdef{\dout^i (t)}{t\!\in
    \!\realnonnegative}, \; i\in{\VV} ,
\end{align*}
where $\lL_{t}$ is the Laplacian of $\GG(t)$. Note that if $\GG$ has
uniformly bounded weights and is weight-balanced and strongly
connected, then $\UHlambda_2>0$.  If $\GG$ is a time-varying connected
graph, we use the notation $\Ulambda_2$ instead of~$\UHlambda_2$.  The
following result, taken from~\citep[Lemma 4.5]{SSK-JC-SM:14-ijrnc}, is
useful when dealing with recurrently jointly strongly connected
digraphs.  When such digraphs are weight-balanced with uniformly
bounded weights,
% For a weight-balanced and recurrently jointly strongly connected
% time-varying digraph with uniformly bounded weights,
there exist $\Hlambda_{\sigma}>0$ and $\rho>0$ such that
\begin{equation}\label{eq::Hlambda_sigma_Def}
  \big \| \text{e}^{-\beta
    {\R}^\top\vect{\mathsf{L}}_{t}{\R}(t-t_0)} \big \|\leq
  \rho \, \text{e}^{-\beta\Hlambda_{\sigma}(t-t_0)},~~~\forall t\geq t_0\geq
  0,
\end{equation}
for any~$\beta>0$.  If the digraph is additionally strongly connected,
then~\eqref{eq::Hlambda_sigma_Def} is satisfied with $\rho=1$ and
$\Hlambda_{\sigma}=\UHlambda_{2}$.

%%%%%%%%%%%%%%%%%%%%%%%%%%%%%%%%%
%%%%%%%%%%%%%%%%%%%%%%%%%%%%%%%%%

\section{Network model and problem
  statement}\label{se:problem-statement}

% This section formalizes the problem of interest.
Consider a network of $N$ agents with single-integrator dynamics,
$\dot{x}^i = g^i$, $i\in\VV$,
% \begin{equation*}%\label{eq::AgentSingInt}
%  \dot{x}^i = g^i,\quad i  \until{N},
%\end{equation*}
where $x^i\in\real$ is the \emph{agreement state} and $g^i\in\real$ is
the \emph{driving command} of agent~$i$.  Our consideration of simple
dynamics is motivated by the fact that the state of the agents does
not necessarily correspond to some physical quantity, but instead to
some logical variable on which agents perform computation and
processing.  Each agent $i \in\VV$ has access to a time-varying reference
signal $\mathsf{r}^i:\realnonnegative\to\real$. Agents transmit
information to other agents through wireless communication and their
interaction topology is modeled by a time-varying weighted
digraph~$\GG$.  An edge $(i,j)$ from $i$ to $j$ at time~$t$ means that
agent $j$ \emph{can} send information to agent $i$ at~$t$.
% Given that communication occurs at discrete instants of time, we let
% $\hat{x}^i$ denote the last known state of agent $i \in\VV$
% transmitted to its in-neighbors.  We let $\{\bar{t}^i_{k}\}_{k \in
% \integersnonnegative}= \{t^i_{k}\}_{k \in \integersnonnegative}\cup
% \{s_k\}_{k \in \integersnonnegative} \subset \realnonnegative$
% denote the sequence of times at which agent $i$ communicates with
% its in-neighbors. Here, $\{t^i_{k}\}_{k \in \integersnonnegative}$
% is the set of times a state to be transmitted is selected, that is
% $\hat{x}^i (t) = x^i(t^i_{k})$ for
% $t\!\in\![{t}^i_{k},{t}^i_{k+1})$, while $\{s_k\}_{k \in
% \integersnonnegative}$ is the set of switching times of the
% underlying communication digraph.  The variable
% $\tilde{x}^i(t)=\hat{x}^i(t)-x^i(t)$ denotes the mismatch between
% the last transmitted state and the state of agent~$i$ at time~$t$.
For convenience, we let $\hat{x}^i$ denote the last sampled state of
agent $i \in\VV$.  We let
\begin{align*}
  \{\bar{t}^i_{k}\}_{k \in \integersnonnegative}= \{t^i_{k}\}_{k \in
    \integersnonnegative}\cup \{s^i_k\}_{k \in \integersnonnegative}
  \subset \realnonnegative
\end{align*}
denote the sequence of \emph{update times} of agent $i \in \VV$.
Here, $\{t^i_{k}\}_{k \in \integersnonnegative}$ is the set of times
at which the state of the agent is sampled, that is $\hat{x}^i (t) =
x^i(t^i_{k})$ for $t\!\in\![{t}^i_{k},{t}^i_{k+1})$, while
$\{s^i_k\}_{k \in \integersnonnegative}$ is the set of switching times
of the underlying communication digraph at which the agent acquires an
in-neighbor. According to this model, agent $i$ broadcasts $\hat{x}^i$
to its in-neighbors at each time $t^i_k$ if $\din^i(t^i_{k}) > 0$ and
also at the times $s^i_k$.  The variable
$\tilde{x}^i(t)=\hat{x}^i(t)-x^i(t)$ denotes the mismatch between the
last sampled state and the state of agent~$i$ at time~$t$.
% We define $\hat{x}^i$ to be the last sampled state of agent
% $i\in\VV$
% taken at time $t^i_k\leq t$, i.e., $\hat{x}^i (t) = x^i(t^i_{k})$
% for
% $t\!\in\![{t}^i_{k},{t}^i_{k+1})$.
  % \begin{definition}[Broadcast policy]\label{def::broadcast}
  %   Agent $i\in\VV$ broadcasts $\hat{x}^i$ at $t={t}^i_{k}$ if
  %   $\din^i(t^i_{k})\neq0$ as well as at any $t\in
  %   \{\{s^i_{k'}\}_{k'
  %   \in \integersnonnegative}\cap ({t}^i_{k},{t}^i_{k+1})\}$.
  % \end{definition}
  % As such, the broadcast times of every agent $i\in\VV$ is a subset of
  % $\{\bar{t}^i_{k}\}_{k \in \integersnonnegative}= \{t^i_{k}\}_{k \in
  %   \integersnonnegative}\cup \{s^i_k\}_{k \in \integersnonnegative}
  % \subset \realnonnegative$. We refer to $\{\bar{t}^i_{k}\}_{k \in
  %   \integersnonnegative}$ as the \emph{set of update times} of agent
  % $i\in\VV$.  The variable $\tilde{x}^i(t)=\hat{x}^i(t)-x^i(t)$
  % denotes the mismatch between the last transmitted state and the
  % state of agent~$i$ at time~$t$.
When the graph is not fixed, the model above assumes that individual
agents are made aware of the identity of newly acquired in-neighbors
(so that they can communicate to them their last transmitted state)
and departed out-neighbors (so that they can remove their state from
their computations). 
% Although we do not get into the specific details of how this service
% is implemented, we envision that a simple protocol based on basic
% messaging can capture agent arrival and departures.

Under the network model described above, our goal is to design a
distributed algorithm that allows each agent to asymptotically track
the average of the reference inputs $\avrg{\mathsf{r}}(t)$ across the
group.  The algorithm design amounts to specifying, for each agent
$i\in\VV$, a suitable distributed driving command
$\map{g^i}{\real^{\mathcal{N}^i}}{\real}$ together with a mechanism
for triggering communication with its in-neighbors in an opportunistic
fashion.  By \emph{distributed}, we mean that each agent only needs to
receive information from its out-neighbors to evaluate~$g^i$ and the
communication triggering law.
% By \emph{opportunistic}, we mean that the transmission of
% information to its in-neighbors should happen at times when it is
% needed to preserve the stability and convergence of the coordination
% algorithm.
A key requirement on the communication triggering mechanism is that
the resulting network evolution is free from Zeno behavior, i.e., does
not exhibit an infinite amount of communication rounds in any finite
amount of time.

\section{Continuous-time computation with distributed event-triggered
  communication}\label{sec::Discrete}
  
Here, we present our solution to the problem stated in
Section~\ref{se:problem-statement}. Our starting point is the
continuous-time algorithm for dynamic average consensus proposed in
our previous work~\citep{SSK-JC-SM:14-ijrnc},
\begin{subequations}\label{eq::CT}
  \begin{align}
    \dot{v}^i & = \alpha\beta\sum\nolimits_{j=1}^N \mathsf{a}_{ij}(x^i-x^j),
    \\
    \dot{x}^i & =
    \!\dot{\mathsf{r}}^i\!-\!\alpha(x^i-\mathsf{r}^i)\!-\!\beta\sum\nolimits_{j=1}^N
    \mathsf{a}_{ij}(x^i-x^j)\!-\!v^i , 
  \end{align}
\end{subequations}
for each $i\in{\VV}$. In~\eqref{eq::CT}, $\alpha, \beta \in
\realpositive$ are design parameters.  Note that the execution of this
algorithm requires continuous agent-to-agent sharing of the
variable~$x$.
% Note that the execution of this algorithm requires agents to
% continuously interchange information about their respective
% variable~$x$ with their neighbors.

\begin{remark}[Knowledge of derivative of reference signals]
  {\rm Interestingly,~\eqref{eq::CT} can be executed without knowledge
    of the time derivatives of the reference signals. In fact, with
    the change of variables $\Bvect{x}=\vect{x}-\vect{\mathsf{r}}$,~\eqref{eq::CT}
    reads
    \begin{align*}
      \dot{v}^i & = \alpha\beta\sum\nolimits_{j=1}^N
      \mathsf{a}_{ij}(\bar{x}^i+ \mathsf{r}^i -\bar{x}^j
      -\mathsf{r}^j),
      \\
      \dot{\bar{x}}^i & = -\alpha
      \bar{x}^i\!-\!\beta\sum\nolimits_{j=1}^N
      \mathsf{a}_{ij}(\bar{x}^i+\mathsf{r}^i-\bar{x}^j
      -\mathsf{r}^j)\!-\!v^i ,
      \\
      x^i&=\bar{x}^i+\mathsf{r}^i,\quad i\in\VV .
    \end{align*}
    Agents operate on their corresponding components of $v$ and
    $\bar{x}$, and interchange with their neighbors the corresponding
    components of~$x$. We only use the representation~\eqref{eq::CT}
    for convenience in our technical analysis later.  } \oprocend
\end{remark}

The following result summarizes, for reference, the asymptotic
correctness guarantees of~\eqref{eq::CT}.

\begin{theorem}[Convergence of~\eqref{eq::CT} over weight-balanced and
  recurrently jointly strongly connected
  digraphs~\citep{SSK-JC-SM:14-ijrnc}]\label{thm::Alg_ContiTime}
  Assume the agent inputs satisfy
  $\|\vect{\Pi}_N\dvect{\mathsf{r}}\|_{\text{ess}} = \gamma<\infty$.
  Let the communication topology be a weight-balanced and recurrently
  jointly strongly connected time-varying digraph~$\GG$ with uniformly
  bounded weights.  Then, for any $\alpha, \beta\! \in\!
  \realpositive$, the evolution of~the algorithm~\eqref{eq::CT}
  over~$\GG$ initialized at $z^i(0),v^i(0)\in \real$ with
  $\SUM{i=1}{N}v^i(0)=0$ is bounded and satisfies, for $i\in\VV$,
  \begin{equation}\label{eq::Alg_D_ultimate_bound} 
    \limsup_{t\to\infty} \Big| x^i(t)-\avrg{\mathsf{r}}(t) \Big| \leq \rho
    \frac{\gamma}{ 
      \beta\Hlambda_{\sigma}},
  \end{equation}
  with $\Hlambda_{\sigma}$ and $\rho$
  satisfying~\eqref{eq::Hlambda_sigma_Def}. 
  % \hl{The rate of convergence is exponential with the
  % least value of $\Hlambda_{\sigma}$}.
\end{theorem}

Given the network model of Section~\ref{se:problem-statement}, where
the transmission of information is limited to discrete instants of
time, we propose here the following implementation of~\eqref{eq::CT}
with discrete-time communication,
\begin{subequations}\label{eq::EventTrig_Alg}
  \begin{align}
    \dot{v}^i & =\alpha\beta\sum\nolimits_{j=1}^N
    \mathsf{a}_{ij}(\hat{x}^i-\hat{x}^j), \label{eq::EventTrig_Alg-a}
    \\
    \dot{x}^i &
    =\dot{\mathsf{r}}^i\!-\!\alpha(x^i-\mathsf{r}^i)\!-\!\beta\sum\nolimits_{j=1}^N 
    \mathsf{a}_{ij}(\hat{x}^i\!-\!\hat{x}^j)\!-\!v^i,
   \label{eq::EventTrig_Alg-b}
  \end{align}
\end{subequations}
for each $i\in{\VV}$.  Our task is to provide individual agents with
triggers that allow them to determine in an opportunistic fashion when
to transmit information to their in-neighbors. The design of such
triggers is challenging because triggers need to be distributed, so
that agents can check them with the information available to them from
their out-neighbors, they must guarantee the absence of Zeno behavior,
and they have to ensure the network achieves dynamic average consensus
even though agents operate with outdated information while inputs are
changing with~time.

\subsection{Compact-form algorithm
  representations}\label{sec::compactform}

Here, we present two equivalent compact-form representations of the
algorithm~\eqref{eq::EventTrig_Alg} for analysis purposes.  First,
consider the change of variable
\begin{subequations}\label{eq::xvToyw}
  \begin{align}
    \vect{w} \! = \!
    \vect{v}-\alpha\vect{\Pi}_N\vect{\mathsf{r}},~~\vect{y}\! =\! \vect{x}-\Bvect{\mathsf{r}} ,~~\Bvect{\mathsf{r}}\!
    = \!\avrg{\mathsf{r}} \vect{1}_N.
    % \label{eq::wTOv}
  \end{align}
\end{subequations}
which transforms the algorithm~\eqref{eq::CT} into (compact form) 
\begin{subequations}\label{eq::Alg_shf}
  \begin{align}
   \dvect{w} &= \alpha\beta\vect{L}_t\vect{y}+
    \alpha\beta\vect{L}_t\Tvect{x}-\alpha\vect{\Pi}_N\dvect{\mathsf{r}},
    \label{eq::Alg_shf-b}\\
    \dvect{y} & = -\alpha\vect{y} - \beta\vect{L}_t \vect{y} -
    \beta\vect{L}_t \Tvect{x}+\vect{\Pi}_N\dvect{\mathsf{r}}-\vect{w},
    \label{eq::Alg_shf-a}
  \end{align}
\end{subequations}
Note that $\vect{y}$ is the aggregated tracking error vector. Here, we have used $\vect{L}_t\Hvect{x} =\vect{L}_t(\vect{x}+\Tvect{x})
=\vect{L}_t\vect{y}+\vect{L}_t\Tvect{x}$ with $\Tvect{x}=\Hvect{x}-\vect{x}$.  
%For the second representation, we introduce $\vect{\mathfrak{T}} \in \real^{N \times N}$, ${\R} \in \real^{N
 % \times N-1}$, and $\vect{r} \in \real^N$ such that
%\begin{equation*}
%  \vect{\mathfrak{T}} 
 % = [
 % \vect{r} , \;{\R}
 % ]
 % ,~~
  %\vect{r} =
  %\frac{1}{\sqrt{N}}\vect{1}_N,~~\vect{r}^\top{\R} =
 % \vect{0},~~{\R}^\top{\R}=\vect{I}_{N-1}.  
%\end{equation*}
%Notice that $\vect{\mathfrak{T}}^\top\vect{\mathfrak{T}}=\vect{I}_N$. 
We use the orthonormal transfer matrix $\vect{\mathfrak{T}} \in \real^{N \times N}$
\begin{equation}\label{eq::T}
  \vect{\mathfrak{T}} 
  \!= \![
  \rr  ~~\R
  ]
  ,~
  \rr =
  \frac{1}{\sqrt{N}}\vect{1}_N,~\vect{r}^\top{\R}\!=\!
  \vect{0},~{\R}^\top{\R}=\vect{I}_{N-1}.  
\end{equation}
to obtain our second representation below that separates out the constant dynamics of the algorithm,
%For the second representation, consider the
%following change of variables, \hl{which is intended to separate the constant dynamics of the algorithm from the rest of the dynamics,}
\begin{equation}\label{eq::trans_stable}
  q_1=\rr^\top\vect{w}, \quad\vect{q}_{2:N}=\alpha
  {\R}^\top \vect{y} +{\R}^\top\vect{w},\quad
  \vect{z}=\vect{\mathfrak{T}}^\top\vect{y}  .
\end{equation}
We partition the new variable $\vect{z}$ as $(z_1, \vect{z}_{2:N})$,
where $z_1\in\reals$. Then, if the network interaction topology is
weight-balanced, the algorithm~\eqref{eq::Alg_shf} can be written as,
\begin{subequations}\label{eq::DEvent_Alg_Separated}
  \begin{align}
    \!\dot{q}_{1} & \!=0,\label{eq::DEvent_Alg_Separated-a}
    \\
    \!\dvect{q}_{2:N} &
    \!=-\alpha\vect{q}_{2:N} \label{eq::DEvent_Alg_Separated-b},
    \\
    \!\dot{z}_{1} & \!=- \alpha z_1-q_1,\label{eq::DEvent_Alg_Separated-c}
    % \\
    % &\dvect{z}_{2:N} =\! -\beta{\R}^\top
    % \vect{L}{\R}\vect{z}_{2:N}\!-\!\beta{\R}^\top\vect{L}{\R}
    % \Tvect{z}_{2:N}\!
    % -\!\vect{q}_{2:N}\!+\!{\R}^\top\dvect{\mathsf{r}}. \label{eq::DEvent_Alg_Separated-d}
    \\
    \!\dvect{z}_{2:N} & \!=\! -\beta{\R}^\top \vect{L}_t{\R}
    \vect{z}_{2:N} \!-\!\beta{\R}^\top
    \vect{L}_t\Tvect{x}\!+\!{\R}^\top\dvect{\mathsf{r}}
    -\!\vect{q}_{2:N} . \label{eq::DEvent_Alg_Separated-d}
  \end{align}
\end{subequations}
We close this section by describing the relationship between the
initial conditions of the variables for each representation. We invoke these relations in our analysis below. Note that
$ \vect{q}_{2:N} = {\R}^{\top}(\alpha \vect{y}+\vect{w}) =
{\R}^{\top}(\alpha (\vect{x}-\vect{\mathsf{r}})+\vect{v})$.  Then,
given $\vect{x}(0),\vect{v}(0)\in\real^N$ with $\sum\nolimits_{i=1}^Nv^i(0) =
0$, and using $\rr^\top\vect{\Pi}_N=\vect{0}$ and
${\R}{\R}^\top = \vect{\Pi}_N = \vect{\Pi}_N^2 $,
\begin{subequations}\label{eq::q0}
  \begin{align}
    q_1(0) &=\rr^\top\vect{w}(0) = \rr^\top\vect{v}(0)=0,
    \\
   \| \vect{q}_{2:N}(0)\| & =\|\alpha\vect{\Pi}_N
    (\vect{x}(0)-\vect{\mathsf{r}}(0))+\vect{v}(0)\|,
    \\
    z_1(0) &=\rr^\top\vect{y}(0) =
    \rr^\top(\vect{x}(0)-\Bvect{\mathsf{r}}(0)),
    \\
    \| \vect{z}_{2:N}(0) \|& = \| \vect{\Pi}_N
    (\vect{x}(0)-\Bvect{\mathsf{r}}(0)) \|.
  \end{align}  
\end{subequations}
%where we have used $\vect{r}^\top\vect{\Pi}_N=\vect{0}$ and ${\R}{\R}^\top = \vect{\Pi}_N = \vect{\Pi}_N^2 $.

\subsection{Communication triggering law for weight-balanced and
  recurrently jointly strongly connected
  digraphs}\label{se:sc-wb-digraphs}

In this section, for networks with time-varying digraph interactions,
we introduce a distributed event-triggered mechanism that agents can
employ to determine their sequence of communication times.  For each
agent, the execution of this mechanism relies on its local variables.
%and the triggered states received from its out-neighbors. 
 This naturally results in asynchronous schedules of
communication, which poses additional analysis
challenges. Nevertheless, we are able to overcome them in the
following result which states that the closed-loop network execution
is free from Zeno behavior and guaranteed to achieve practical dynamic
average consensus. 

\begin{theorem}[Convergence of~\eqref{eq::EventTrig_Alg} over
  recurrently jointly strongly connected and weight-balanced digraph
  with asynchronous distributed event-triggered
  communication]\label{thm::Alg_D_Event}
  Assume that the input of each agent $i\in\VV$ satisfies
  $|\dot{\mathsf{r}}^i|_{\text{ess}}\!=\! \kappa^i<\infty$, while
  $\|\vect{\Pi}_N\dvect{\mathsf{r}}\|_{\text{ess}} = \gamma<\infty$.
  Let the communication topology be a weight-balanced and recurrently
  jointly strongly connected time-varying digraph~$\{\GG\}_{\{s_k\}}$
  with uniformly bounded weights.  For $\vect{\eps} \in
  \realpositive^N$, consider an implementation of the
  algorithm~\eqref{eq::EventTrig_Alg} over $\{\GG\}_{\{s_k\}}$, where
  agents communicate according to the model described in
  Section~\ref{se:problem-statement} with the sampling times
  $\{t^i_{k}\}_{k \in \integersnonnegative}$ of agent $i\!\in\!\VV$,
  starting at
  % follows the broadcast policy given in
  % Definition~\ref{def::broadcast} with
  % communicates with its neighbors at times
  % $\{\bar{t}^i_{k}\}_{k \in \integersnonnegative}= \{t^i_{k}\}_{k
  % \in
  % \integersnonnegative} \cup \{s^i_k\}_{k \in \integersnonnegative}
  % \subset \realnonnegative$, where $\{s_k\}_{k}$ is the set of
  % switching times of the digraph and
  % $\{t^i_{k}\}$, started at
  $t^i_0=0$, determined~by
  \begin{align}\label{eq::TrigLaw_Distributed_Directed}
    t^i_{k+1} \!=\! \argmax \setdef{t \in [t^i_{k},\infty)\!}{\!
      |{x}^i(t^i_{k})\!-\!x^i(t)|\leq \eps_i }.
  \end{align}
  % Recall $\{s_k\}_{k \in \integersnonnegative}$ is the set of
  % switching times of the digraph.
  Then, for any $\alpha, \beta\! >0$, the evolution starting from
  $x^i(0)\in\reals$ and $v^i(0)\in\reals$ with
  $\sum\nolimits_{i=1}^N\!v^i(0)\!=\!0$ satisfies
  \begin{align}\label{eq::Alg_D_Event_ultimate_bound}
    \limsup_ {t\to\infty}\! \Big|x^i(t)\!-\!\avrg{\mathsf{r}}(t) \Big|
    \!\leq \!\frac{\gamma\!+\!\beta
      \BnL\,\|\vect{\eps}\|}{\beta{\Hlambda}_{\sigma}}\rho ,
  \end{align}
  for $ i\in\VV$ with an exponential rate of convergence of
  $\min\{\alpha,\beta{\Hlambda}_{\sigma}\}$. Here, $\Hlambda_{\sigma}$
  and $\rho$ satisfy~\eqref{eq::Hlambda_sigma_Def}. Furthermore, the
  inter-execution times of~\eqref{eq::TrigLaw_Distributed_Directed}
  for each agent $i\in\VV$ are lower bounded by
  \begin{equation}\label{eq::dist_event_time_lw_bnd}
    \tau^i = \frac{1}{\alpha} \ln \Big(1+\frac{\alpha
      \eps_i}{c^i} \Big) ,
  \end{equation}  
  where
  \begin{align}\label{eq::ci}
    c^i &= \kappa^i+(\alpha + 2\beta \Bdout^i)
    \sqrt{\eta^2+|\vect{r}^\top(\vect{x}(0)-\Bvect{\mathsf{r}}(0))|^2}
    \nonumber
    \\
    & \quad + \|\vect{\Pi}_N(\alpha
    (\vect{x}(0)-\vect{\mathsf{r}}(0))+\vect{v}(0))\| + \alpha\eta,
  \end{align}
  and
  \begin{align*}
    \left.\begin{array}{ll} \eta & =
        \rho\frac{(\gamma+\beta\BnL\,\|\vect{\eps}\|^2)}{\beta{\Hlambda}_{\sigma}}\!
        +\!  \|\vect{\Pi}_N (\vect{x}(0)-\Bvect{\mathsf{r}}(0))\|\!+\!
        \rho\|\vect{q}_{2:N}(0) \| \times
        \\
        & \quad
    \begin{cases}
      \frac{1}{\alpha-\beta{\Hlambda}_{\sigma}}((
      \frac{\beta{\Hlambda}_{\sigma}}{
        \alpha})^{\frac{\beta{\Hlambda}_{\sigma}}{\alpha-\beta{\Hlambda}_{\sigma}}}
      - (\frac{\beta{\Hlambda}_{\sigma}}{
        \alpha})^{\frac{\alpha}{\alpha-\beta{\Hlambda}_{\sigma}}}), &
      \text{if }\beta{\Hlambda}_{\sigma}\neq \alpha ,
      \\
      \frac{1}{\beta{\Hlambda}_{\sigma} \text{e}} , & \text{if
      }\beta{\Hlambda}_{\sigma}= \alpha .
    \end{cases}
    \end{array}\right.
  \end{align*}
  Hence, $\{\bar{t}_k\}_{k \in \integersnonnegative} = \cup_{i=1}^N
  \cup_{k \in \integersnonnegative} \bar{t}^i_{k}$
  % $\{\bar{t}^i_{k}\}_{k \in \integersnonnegative}$ , the set of
  % times in the dynamics of each agent $i\in\VV$
  has no accumulation point and the execution
  of~\eqref{eq::EventTrig_Alg} over $\{\GG\}_{\{s_k\}}$ is Zeno-free.
  % (recall $\{\bar{t}^i_{k}\}_{k \in \integersnonnegative}=
  % \{t^i_{k}\}_{k \in \integersnonnegative}\cup \{s^i_k\}_{k \in
  % \integersnonnegative}$).
\end{theorem}
\begin{pf}
  Given an initial condition, let $[0,T)$ be the maximal interval on
  which there is no accumulation point in the set of update times
  $\{\bar{t}_k\}_{k \in \integersnonnegative} $. 
  %= \cup_{i=1}^N \cup_{k\in \integersnonnegative} \bar{t}^i_{k} $.  
  Note that $T>0$, since
  the number of agents is~finite and, for each $i\in\VV$, $\eps_i>0$
  and $\tilde{x}^i(0) = \hat{x}^i(0)-x^i(0)=0$.  The
  dynamics~\eqref{eq::EventTrig_Alg}, under the event-triggered
  communication scheme~\eqref{eq::TrigLaw_Distributed_Directed}, has a
  unique solution in the time interval $[0,T)$.  Our first step is to
  show that the trajectory stays bounded in $[0,T)$.  Consider the
  compact-form representation~\eqref{eq::DEvent_Alg_Separated} of the
  algorithm.  Given $\sum\nolimits_{i=1}^Nv^i(0)=0$
  and~\eqref{eq::q0}, for $t\in\realnonnegative$, \smallskip
  \begin{align}\label{eq::a_to_c_Solution}
   \!\! q_1(t) \!=\! 0, \; \vect{q}_{2:N}(t) \!=\!
    \vect{q}_{2:N}(0)\text{e}^{-\alpha t}\!\!, \;
    z_1(t) \!=\!z_1(0)\text{e}^{-\alpha t}\!.
  \end{align}
  Hence these variables are bounded. To bound $t \mapsto
  \vect{z}_{2:N}(t)$, we look into the solution
  of~\eqref{eq::DEvent_Alg_Separated-d} by substituting
  $\vect{q}_{2:N}(t) = \vect{q}_{2:N}(0)\text{e}^{-\alpha t}$ and
  considering $(\Tvect{x},{\R}\dvect{\mathsf{r}})$ as exogenous
  inputs. In the time interval $t\in[0,T)$ that this solution exists,
  one has
  \begin{align*}
    \vect{z}_{2:N}(t) &=
    \vect{\Phi}(t,0)\vect{z}_{2:N}(0)+\int_{0}^t\vect{\Phi}(t,\tau)
    \e^{-\alpha\tau}\vect{q}_{2:N}(0) \text{d}\tau \nonumber
    \\
    &-\int_{0}^t\vect{\Phi}(t,\tau)(\beta{\R}^\top
    \vect{L}_{\tau}\Tvect{x}(\tau)-{\R}^\top\dvect{\mathsf{r}}(\tau))
    \text{d}\tau
    % {\R}\vect{L}_{\sigma(\tau)}\vect \text{d}\tau
  \end{align*}
  where $\vect{\Phi}(t,\tau)=\e^{-\beta{\R}^\top
    \vect{L}_{t}{\R}(t-\tau)}$.  Given the
  event-triggered communication
  law~\eqref{eq::TrigLaw_Distributed_Directed}, we have
  $\|\Tvect{x}\|\leq \|\vect{\eps}\|$. Then, for $t\in[0,T)$ we have
  \begin{align}\label{eq::z2N_bound_exact}
    \|\vect{z}_{2:N}(t)\|\leq&
    \rho\text{e}^{-\beta\Hlambda_{\sigma}t}\|\vect{z}_{2:N}(0)\|+
    \\
    &\rho \frac{\gamma + \beta\BnL\,
      \|\vect{\eps}\|}{\beta{\Hlambda}_{\sigma}}(1-
    \text{e}^{-\beta{\Hlambda}_{\sigma} t})+
    \rho\|\vect{q}_{2:N}(0)\|\times \nonumber
    \\
    & \begin{cases} \frac{ 1}{\alpha-\beta{\Hlambda}_{\sigma}}
      (\text{e}^{-\beta{\Hlambda}_{\sigma} t}-\text{e}^{-\alpha t}), &
      \text{if }\beta{\Hlambda}_{\sigma}\neq \alpha ,
      \\
      t \text{e}^{-\beta\Hlambda_\sigma t}, & \text{if
      }\beta{\Hlambda}_{\sigma}= \alpha .
    \end{cases}\nonumber
  \end{align}
  Here, we also used
  $\|{\R}^\top\dvect{\mathsf{r}}\|=\|\vect{\Pi}\dvect{\mathsf{r}}\|\leq
  \gamma$ and~\eqref{eq::Hlambda_sigma_Def}.  Taking the maximum of
  each term in the righthand side of~\eqref{eq::z2N_bound_exact} and
  using~\eqref{eq::q0}, we have
  \begin{align}\label{eq::z2N_bound_D}
    \|\vect{z}_{2:N}(t)\|\leq \eta , \quad t\in[0,T) ,
  \end{align}
  were the constant $\eta$ is given in the statement.
  %
  % \begin{align}\label{eq::z2N_bound_D}
  %   & \|\vect{z}_{2:N}(t)\|\leq \eta
  %   \frac{(\gamma+\beta\BnL \|\vect{\eps}\|^2)}{\beta\UHlambda_2} +
  %   \|\vect{z}_{2:N}(0) \|+
  %   \\
  %   & \begin{cases} \frac{1}{\alpha-\beta\UHlambda_2}
  %     \|\vect{q}_{2:N}(0)
  %     \|((\frac{\beta\UHlambda_2}{\alpha})^{
  %       \frac{\beta\UHlambda_2}{\alpha-\beta\UHlambda_2}}
  %     -
  %     (\frac{\beta\UHlambda_2}{\alpha})^{
  %       \frac{\alpha}{\alpha-\beta\UHlambda_2}})
  %     & \text{if }\beta\UHlambda_2\neq \alpha
  %     \\
  %     \frac{\text{e}^{-1}}{\beta\UHlambda_2}\|\vect{q}_{2:N}(0)\| &
  %     \text{if
  %     }\beta\UHlambda_2= \alpha
  %   \end{cases}\nonumber
  % \end{align}
  % Recall~\eqref{eq::q0}, then we can conclude that the right-hand
  % side of~\eqref{eq::z2N_bound_D} is equal to
  % the constant $\eta$ given in the statement.

  Given that the number of the agents is finite, $\{\bar{t}_k\}_{k \in
    \integersnonnegative}$ is free of Zeno if both $\{s_k\}_{k \in
    \integersnonnegative}$ and $\{t^i_k\}_{k \in
    \integersnonnegative}$, $i \until{n}$, are uniformly lower
  bounded. Since the former fact holds by assumption, we next
  establish the latter. That is, for each agent, we establish a lower
  bound on the inter-execution times
  of~\eqref{eq::TrigLaw_Distributed_Directed} by determining a lower
  bound on the amount of time it takes for agent $i\in\VV$ to have
  $|\hat{x}^i-x^i|$ evolve from $0$ to $\eps_i$. Note that
  \begin{align}\label{eq:bound-y}
    \|\vect{y}(t) \|= \|\vect{z}
    (t)\| % = \sqrt{\|\vect{z}_{2:N}\|^2+|z_1|^2}
    \le
    \sqrt{\eta^2\!+\!|\vect{r}^\top(\vect{x}(0)-\Bvect{\mathsf{r}}(0))|^2}
    ,
  \end{align}
  for all $t \in [0,T)$, where we have
  used~\eqref{eq::q0},~\eqref{eq::a_to_c_Solution}
  and~\eqref{eq::z2N_bound_D}.  On the other hand, since
  $\vect{r}^\top\vect{w}(t)=0$ for all $t\in\realnonnegative$
  by~\eqref{eq::a_to_c_Solution}, we have $\pPi_N \vect{w}(t) =
  \vect{w}(t)$. Multiplying the second equation
  in~\eqref{eq::trans_stable} by ${\R}$ and using ${\R}{\R}^\top =
  \pPi_N$, we obtain $\vect{w}(t)={\R}\vect{q}_{2:N}(t) -
  \alpha{\R}\vect{z}_{2:N}(t)$. Using now
  ${\R}^\top{\R}=\vect{I}_{N-1}$, we deduce
  \begin{align*}
    \|\vect{w} (t) \| & \leq
    \|\vect{q}_{2:N}(t)\|+\alpha\|\vect{z}_{2:N}(t) \|
    \\
    & \le \|\alpha\vect{\Pi}_N(\vect{x}(0)-\vect{\mathsf{r}}(0)) +
    \vect{v}(0)\|+\alpha\eta,
  \end{align*}
  for all $t\in[0,T)$, where we have again
  used~\eqref{eq::q0},~\eqref{eq::a_to_c_Solution}
  and~\eqref{eq::z2N_bound_D}. Next notice that, for each $i\in\VV$,
  using~\eqref{eq::EventTrig_Alg} and~\eqref{eq::xvToyw}, we have
  \begin{align*}
    & \frac{d}{dt} |\hat{x}^i-x^i|=
    -\frac{(\hat{x}^i-x^i)^\top\dot{x}^i}{|\hat{x}^i-x^i|} \leq
    |\dot{x}^i|
    \\
    & = |\dot{\mathsf{r}}^i -\alpha
    (x^i-\mathsf{r}^i)-\beta\sum\nolimits_{j=1}^N
    \mathsf{a}_{ij}(\hat{x}^i-\hat{x}^j)-v^i|
    \\
    & = |\dot{\mathsf{r}}^i \!-\!\alpha
    (x^i-\avrg{\mathsf{r}})\!-\!\beta\sum\nolimits_{j=1}^N
    \mathsf{a}_{ij}(\hat{y}^i\!-\!\hat{y}^j)-w^i|
    \\
    & \leq \kappa^i \!+\!
    \alpha|\hat{x}^i-x^i|\!+\!\alpha|\hat{y}^i|\!+\!\beta\sum\nolimits_{j=1}^N
    \mathsf{a}_{ij}(|\hat{y}^i|\!+\!|\hat{y}^j|)\!+\!|w^i| .
    % & \leq
    % \kappa^i+\alpha|\hat{x}^i-x^i| +
    % (\alpha+2{\mathsf{d}}_{\text{out}}^i)\sqrt{\eta^2
    %   +|\vect{r}^\top(\vect{x}(0)
    %   - \Bvect{x}(0))|^2}+|w^i|
  \end{align*}
  Then, by noting that $|y^i|\leq\|\vect{y}\|$ and
  $|w^i|\leq\|\vect{w}\|$, and using the bounds established above on
  $\|\vect{y}\|$ and $\|\vect{w}\|$, we obtain
  \begin{align}\label{eq::dot_x_bound}
    \frac{d}{dt} |\hat{x}^i-x^i| &\leq \alpha|\hat{x}^i-x^i|+c^i,
  \end{align}
  where $c^i$ is given in~\eqref{eq::ci}.
% The application of
% Lemma~\ref{lem::growth_xhat} now yields
% \begin{align}\label{eq::dot_x_bound}
%   \frac{d}{dt} |\hat{x}^i-x^i| &\leq \alpha|\hat{x}^i-x^i|+c^i,
% %   & \leq
% %   \kappa^i+\alpha|\hat{x}^i-x^i|+(\alpha +
   % % 2{\dout^i)\sqrt{\eta^2
   % %   + |\vect{r}^\top(\vect{x}(0)-\Bvect{x}(0))|^2}+|w^i|
  %\end{align}
 % where $c^i$ is given in~\eqref{eq::ci}.
  % $c^i = \big(\kappa^i+(\alpha+2{\dout^i)\sqrt{\eta^2 +
  %     |\vect{r}^\top(\vect{x}(0)-\Bvect{x}(0))|^2} +
  %   \|\alpha\vect{\Pi}_N(\vect{x}(0)-\vect{\mathsf{r}}(0))+\vect{v}(0)\|
  %   + \alpha\eta \big)$.
Next, using the Comparison Lemma, cf.~\cite{HKK:02}, and $\hat{x}^i=
{x}^i(t^i_{k})$, we deduce
  \begin{align}\label{eq::state_holding_error}
    |\hat{x}^i-{x}^i (t)| \leq \frac{c^i}{\alpha}(\text{e}^{\alpha
      (t-t^i_{k})}-1), \quad t \ge t^i_{k} .
  \end{align}
  Therefore, the time it takes $|\hat{x}^i-{x}^i|$ to reach $ \eps_i$
  is lower bounded by $\tau^i>0$ as given
  in~\eqref{eq::dist_event_time_lw_bnd}.  This fact also implies that
  $T=\infty$. To see this, we reason by contradiction, i.e., suppose
  $T<\infty$. Then, the sequence of times
  $\{\bar{t}_k\}_{k\in\integersnonnegative}$ has an accumulation point
  at $T$. Because there is a finite number of agents, this implies
  that there is an agent $i\in{\VV}$ for which
  $\{\bar{t}^i_{k}\}_{k\in\integersnonnegative}$ has an accumulation
  point at~$T$. This implies that $i$ transmits infinitely often in
  the time interval $[T-\Delta,T)$ for any $\Delta\in(0,T]$. Given
  that the switching times of the digraph are uniformly lower bounded,
  this contradicts the fact that the inter-event times
  of~\eqref{eq::TrigLaw_Distributed_Directed} are lower bounded by
  $\tau^i>0$ on $[0,T)$.
  % Having established $T=\infty$, note that this fact implies that
  % under the event-triggered communication
  % law~\eqref{eq::TrigLaw_Distributed_Directed}, the
  % algorithm~\eqref{eq::EventTrig_Alg} does not exhibit Zeno
  % behavior.
  Finally, note that
  \begin{align}\label{eq:auxx}
    |x^i-\avrg{\mathsf{r}}| & \le \|\vect{x}-\Bvect{\mathsf{r}} \| =
    \|\vect{y} \| = \|\vect{z}\|,
  \end{align}
  for $i\in{\VV}$, so the bound~\eqref{eq:bound-y} holds for all $t
  \in \realnonnegative$.  From~\eqref{eq::a_to_c_Solution},
  $\lim_{t\to\infty}z_1(t)=0$ with an exponential convergence rate
  $\alpha$. From~\eqref{eq::z2N_bound_exact},
  $\limsup_{t\to\infty}\|z_{2:N}(t)\| \leq
  {(\gamma+\beta\BnL \|\vect{\eps}\|)}/{(\beta\UHlambda_2)}$ with an
  exponential rate of convergence $\min\{\alpha,\beta\UHlambda_2\}$.
  Combining these facts with~\eqref{eq:auxx}
  yields~\eqref{eq::Alg_D_Event_ultimate_bound}, and this concludes
  the proof.  \boxend
\end{pf}

% In the following, we provide several observations about the execution of
% the algorithm~\eqref{eq::EventTrig_Alg} under the proposed
% event-triggered communication
% law~\eqref{eq::TrigLaw_Distributed_Directed}.  
Not surprisingly, the ultimate convergence error
bound~\eqref{eq::Alg_D_Event_ultimate_bound} 
% obtained under event-triggered discrete-time communication
is worse than the bound~\eqref{eq::Alg_D_ultimate_bound} obtained when
agents communicate continuously.  The
trigger~\eqref{eq::TrigLaw_Distributed_Directed} does not use the full
agent state and hence can be seen as an output feedback
event-triggered controller, see e.g.,~\citep{MCFD-WPMHH:12}, for which
guaranteeing the existence of lower-bounded inter-execution times is
in general difficult.

% The combination of the facts that the total number of agents is finite
% and each agent's inter-event times are lower bounded implies that the
% total number of events in any finite time interval is finite. In
% general, an explicit expression lower bounding the network inter-event
% times is not available.

\begin{remark}[Inter-event times as a function of the design
  parameters]\label{rem:tau-dependence}
  {\rm The lower bound $\tau^i$ in~\eqref{eq::dist_event_time_lw_bnd}
    along with knowledge of the lower bound on the switching times of
    the topology allows to compute bounds on the maximum number of
    communication rounds (and associated energy spent) by each agent
    $i \in\VV$ (and hence the network) during a given time interval.
    This lower bound depends on the various ingredients as follows:
    $\tau^i$ is an increasing function of $\eps_i$ and a decreasing
    function of $\alpha$ and $c^i$. Through the latter variable, the
    bound also depends on the graph topology and the design parameter
    $\beta$.  Given the definition of $c^i$, one can deduce that the
    faster an input of an agent is changing (larger $\kappa^i$) or the
    farther the agent initially starts from the average of the inputs,
    the more often that agent would need to trigger communication.
    The connection between network performance and communication
    overhead can also be noted here. Increasing $\beta$ or decreasing
    $\eps_i$ to improve the
    bound~\eqref{eq::Alg_D_Event_ultimate_bound} results in smaller
    inter-event times.  Given that the convergence rate
    of~\eqref{eq::EventTrig_Alg}
    under~\eqref{eq::TrigLaw_Distributed_Directed}
    is~$\min\{\alpha,\beta\Hlambda_{\sigma}\}$, decreasing $\alpha$ to
    increase the inter-event times slows down convergence.  \oprocend
  }
\end{remark}

\subsection{Communication triggering law for time-varying connected
  undirected graphs}\label{se:c-graphs}
 
Here, we design a distributed event-triggered communication law
for~\eqref{eq::EventTrig_Alg} over networks with time-varying
connected undirected graph interaction topologies.  While the results
of the previous section are valid for these topologies, here we show
that the structural properties of the Laplacian matrix in the
undirected case allows the alternative event-triggered law to have
longer inter-event times with similar tracking performance.

\begin{proposition}[Convergence of~\eqref{eq::EventTrig_Alg} over
  time-varying connected undirected graphs with asynchronous
  distributed event-triggered
  communication]\label{prop::Alg_D_Event_UD}
  Assume that the input of each agent $i\in{\VV}$ satisfies
  $|\dot{\mathsf{r}}^i|_{\text{ess}}= \kappa^i<\infty$, while
  $\|\vect{\Pi}_N\dvect{\mathsf{r}}\|_{\text{ess}} = \gamma<\infty$.
  Let the communication topology be a connected, piecewise continuous
  time-varying undirected graph~$\{\GG\}_{s_k}$ with uniformly bounded
  weights.  For $\vect{\eps} \in \realpositive^N$, consider an
  implementation of the algorithm~\eqref{eq::EventTrig_Alg} over
  $\{\GG\}_{\{s_k\}}$, where agents communicate according to the model
  described in Section~\ref{se:problem-statement} with the sampling
  times $\{t^i_{k}\}_{k \in \integersnonnegative}$ of agent
  $i\!\in\!\VV$, starting at $t^i_0=0$, determined~by
  \begin{align}\label{eq::TrigLaw_Distributed_UD}
    t^i_{k+1} &= \argmax \setdef{t \in [t^i_{k},\infty)}{|\hat{x}^i(t)
      -x^i(t) |^2
      \\
      & \leq \frac{1}{4 \dout^i(t) } \sum\nolimits_{j=1}^N
      {\mathsf{a}}_{ij}(t)|\hat{x}^i(t) - \hat{x}^j(t)|^2 \!+\!
      \frac{\eps_i^2}{4 \dout^i(t) } }. \nonumber
  \end{align}
  Then, for any $\alpha, \beta \in \realpositive$, the evolution
  starting from $x^i(0)\in\reals$ and $v^i(0)\in\reals$ with
  $\sum\nolimits_{i=1}^Nv^i(0)=0$ satisfies
  \begin{align}\label{eq::Alg_D_Event_ultimate_bound_UD}
    \limsup_ {t\to\infty} \Big| x^i(t)-\avrg{\mathsf{r}}(t) \Big| &
    \leq \frac{\gamma}{\beta\Ulambda_2} +
    \\
    & \quad \sqrt{\Big(\frac{\gamma}{ \beta\Ulambda_2}\Big)^2 +
      \frac{\|\vect{\eps}\|^2}{2\Ulambda_2}}, \nonumber
  \end{align}
  % \begin{align}\label{eq::Alg_D_Event_ultimate_bound_UD}
   % \limsup_ {t\to\infty} \Big| x^i(t)-\avrg{\mathsf{r}}(t) \Big| &
   % \leq \frac{\gamma}{\beta\Ulambda_2\theta} +
   % \\
  %  & \quad \sqrt{\Big(\frac{\gamma}{ \beta\Ulambda_2\theta}\Big)^2 +
   %   \frac{\|\vect{\eps}\|^2}{2\Ulambda_2\theta}}, \nonumber
  %\end{align}
  %
  %
  for $ i\in\VV$. % and some $0<\theta<1$.
  Furthermore, the inter-execution times of agent $i\in{\VV}$ are
  lower bounded by
  \begin{equation}\label{eq::dist_event_time_lw_bnd_undirected}
    \tau^i = \frac{1}{\alpha} \ln \Big(1+\frac{\alpha
      \eps_i}{2c^i\sqrt{\Bdout^i}} \Big) ,
  \end{equation}  
  where $c^i$ is given in~\eqref{eq::ci},  with $\eta$ substituted by
  \begin{equation*}
    \left.
      \begin{array}{ll}
        \zeta =& \max\Big\{
        \|\vect{\Pi}_N(\vect{x}(0) -\Bvect{\mathsf{r}}(0))\|,
        \frac{\alpha\|\vect{\Pi}_N (\vect{x}(0)-\vect{\mathsf{r}}(0)) +
          \vect{v}(0)\|}{2}+ 
        \\ 
        ~& \frac{\gamma}{\beta\Ulambda_2} +
        \sqrt{(\frac{\|\alpha\vect{\Pi}_N
            (\vect{x}(0)-\vect{\mathsf{r}}(0))+\vect{v}(0)\|}{2} +
          \frac{\gamma}{\beta\Ulambda_2})^2 +
          \frac{\|\vect{\eps}\|^2}{2\Ulambda_2}}\Big\}.
      \end{array}
    \right.
  \end{equation*}
  %\begin{equation*}
   % \left.
    %  \begin{array}{ll}
     %   \zeta =& \max\Big\{
     %   \|\vect{\Pi}_N(\vect{x}(0) -\Bvect{\mathsf{r}}(0))\|,
      %  \frac{\alpha\|\vect{\Pi}_N (\vect{x}(0)-\vect{\mathsf{r}}(0)) +
      %    \vect{v}(0)\|}{2\theta}+ 
       % \\ 
     %   ~& \frac{\gamma}{\beta\Ulambda_2\theta} +
      %  \sqrt{(\frac{\|\alpha\vect{\Pi}_N
       %     (\vect{x}(0)-\vect{\mathsf{r}}(0))+\vect{v}(0)\|}{2\theta} +
        %  \frac{\gamma}{\beta\Ulambda_2\theta})^2 +
        %  \frac{\|\vect{\eps}\|^2}{2\Ulambda_2\theta}}\Big\}.
      %\end{array}
    %\right.
  %\end{equation*}
  Hence, $\{\bar{t}_k\}_{k \in \integersnonnegative} = \cup_{i=1}^N
  \cup_{k \in \integersnonnegative} \bar{t}^i_{k}$
  % $\{\bar{t}^i_{k}\}_{k \in \integersnonnegative}$ , the set of
  % times in the dynamics of each agent $i\in\VV$
  has no accumulation point and the execution
  of~\eqref{eq::EventTrig_Alg} over $\{\GG\}_{\{s_k\}}$ is Zeno-free.
  % (recall $\{\bar{t}^i_{k}\}_{k \in \integersnonnegative}=
  % \{t^i_{k}\}_{k \in \integersnonnegative}\cup \{s^i_k\}_{k \in
  %   \integersnonnegative}$).
\end{proposition}
\begin{pf} 
  % The proof follows the similar steps as of the proof of
  % Theorem~\ref{thm::Alg_D_Event}.
  Given an initial condition, let $[0,T)$ be the maximal interval on
  which there is no accumulation point in the set of event times
  $\{\bar{t}_k\}_{k \in \integersnonnegative} = \cup_{i=1}^N \cup_{k
    \in \integersnonnegative} \bar{t}^i_{k} $.  The expressions
  in~\eqref{eq::a_to_c_Solution} are equally valid in this case. To
  bound $t \mapsto \vect{z}_{2:N}(t)$,
  consider% the candidate Lyapunov
  % function
  \begin{align*}%\label{eq::Lyapun-Digraph}
    V(\vect{z}_{2:N}) &= \frac{1}{2}\vect{z}_{2:N}^\top
    \vect{z}_{2:N}.
  \end{align*}
  The derivative of $V(\vect{z}_{2:N}) $ along the trajectories
  of~\eqref{eq::DEvent_Alg_Separated-d} can be upper bounded, for
  $t\in[0,T)$, as
  \begin{align*}
   & \dot{V}
    % &-\vect{z}_{2:N}^\top\vect{q}_{2:N}(0)e^{-\alpha
    %   t}-\beta\vect{z}_{2:N}^\top{\R}^\top\vect{L}{\R}\vect{z}_{2:N}-
    % \\
    % &
    % \beta\vect{z}_{2:N}^\top{\R}^\top\vect{L}{\R}\Tvect{z}_{2:N}
    % + \vect{z}_{2:N}^\top{\R}^\top\dvect{u}
    % \\
   \leq \|\vect{z}_{2:N}\|\|\vect{q}_{2:N}(0)\|e^{-\alpha t}-
   \frac{1}{2}\beta\Ulambda_2\vect{z}_{2:N}^\top\vect{z}_{2:N}-
   \\
   & \frac{1}{2}\beta(\vect{z}_{2:N}^\top{\R}^\top
   \vect{L}_t{\R}\vect{z}_{2:N} +
   2\vect{z}_{2:N}^\top{\R}^\top\vect{L}_t\Tvect{x}) +
   \gamma\|\vect{z}_{2:N}\|.
  \end{align*}
  %
  %
  % The second term in right-hand side of the
  % inequality above is obtained using $\Ulambda_2 \vect{I}\leq
  % {\R}^\top\vect{L}{\R} $.
  For convenience, let
%  \begin{equation*}
$
    s \!=\! - \vect{z}_{2:N}^\top{\R}^\top\! \vect{L}_t{\R}
    \vect{z}_{2:N}  - 2 \vect{z}_{2:N}^\top{\R}^\top\!
    \vect{L}_t\Tvect{x} $.
 % \end{equation*}
  Using ${\R} {\R}^\top=\vect{\Pi}_N$,
  $\vect{L}_t\vect{\Pi}_N=\vect{\Pi}_N\vect{L}_t=\vect{L}_t$,
  and~\eqref{eq::trans_stable}, we obtain
    $s = -\vect{x}^\top\vect{L}_t\vect{x}-2\vect{x}^\top\vect{L}_t\Tvect{x}
   $. Given that the communication graph is connected for all $t$, we can show (see~\cite{SSK-JC-SM:14-ijrnc} for details)
  \begin{align*}
    s & \le \frac{1}{2}\sum\nolimits_{i=1}^N \!\!\Big(
    4\mathsf{d}_{\text{out}}^i(t)|\hat{x}^i\!-\!{x}^i|^2 \!-\!\sum\nolimits_{j=1}^N
    \mathsf{a}_{ij}(t)|\hat{x}^i\!-\!\hat{x}^j|^2\Big),
  \end{align*}
  which, together with~\eqref{eq::TrigLaw_Distributed_UD} implies that $s\leq
   \frac{1}{2}\|\vect{\eps}\|^2$ for $t \in [0,T)$ under the
   event-triggered communication
   law~\eqref{eq::TrigLaw_Distributed_UD}.  Therefore, for $\theta \in
   (0,1)$, we have
  \begin{align}\label{eq::dV_intermidate}
    \dot{V} \!\!& \leq\! \!- \frac{\beta\Ulambda_2(1\!-\!\theta)}{2}
    \vect{z}_{2:N}^\top\vect{z}_{2:N} \!\!+\! \frac{\beta\Ulambda_2}{2} \Big(\!
    \frac{2}{\beta\Ulambda_2}\|\vect{z}_{2:N}\| \|\vect{q}_{2:N}(0)\|e^{\!-\alpha
      t}\nonumber
    \\
    & \quad
    -\theta\vect{z}_{2:N}^\top\vect{z}_{2:N}+\frac{1}{2\Ulambda_2}
    \|\vect{\eps}\|^2+\frac{2\gamma}{\beta\Ulambda_2}
    \|\vect{z}_{2:N}\| \Big) \nonumber
    \\
    & \leq
    -\frac{\beta\Ulambda_2}{2}(1-\theta)\vect{z}_{2:N}^\top\vect{z}_{2:N}+
    \frac{\beta\Ulambda_2}{2}r,
  \end{align}
  where $r= \frac{2}{\beta\Ulambda_2}\|\vect{z}_{2:N}\|\,
  \|\vect{q}_{2:N}(0)\| - \theta\vect{z}_{2:N}^\top\vect{z}_{2:N} +
  \frac{1}{2\Ulambda_2}\|\vect{\eps}\|^2 +
  \frac{2\gamma}{\beta\Ulambda_2}\|\vect{z}_{2:N}\|$.
  % Here we used
  % $\|\vect{q}_{2:N}(0)\|e^{-\alpha t}\leq \|\vect{q}_{2:N}(0)\|$ for
  % all $t\in \realnonnegative$.
  Notice that $r<0$ for
  \begin{align*}%\label{eq::r_sign_neg}
    \|\vect{z}_{2:N}\|\geq &
    \frac{\|\vect{q}_{2:N}(0)\|}{\beta\Ulambda_2\theta} \!+\!
    \frac{\gamma}{\beta\Ulambda_2\theta} \nonumber
    \\
    & + \sqrt{\Big(\frac{\|\vect{q}_{2:N}(0)\|}{\beta\Ulambda_2\theta}
      + \frac{\gamma}{\beta\Ulambda_2\theta}\Big)^2 +
      \frac{\|\vect{\eps}\|^2}{2\Ulambda_2\theta}} = \hat{\zeta}.
  \end{align*}
  Hence, for $t\in[0,T)$, as long as $\|\vect{z}_{2:N} (t) \|\geq
  \hat{\zeta}$, one has
  \begin{align*}
    \dot{V}\leq &-\frac{1}{2}\beta
    \Ulambda_2(1-\theta)\vect{z}_{2:N}^\top\vect{z}_{2:N} .
  \end{align*}
  Combining this inequality with the definition of $V$, and
  considering the limiting case $\theta\to 1$, we deduce that, for any
  $\vect{z}_{2:N}(0)\in\real^{N-1}$ and $t\in[0,T)$,
  \begin{equation}\label{eq::z2N_bound}
    \|\vect{z}_{2:N}(t)\|\leq \max\{ \|\vect{z}_{2:N}(0)\|, \hat{\zeta}\}=\zeta.
  \end{equation}
  Next, following the same arguments as in the proof of
  Theorem~\ref{thm::Alg_D_Event}, one can establish a lower bound on
  the inter-execution times of any agent $i\in\VV$. To do this, we
  determine a lower bound on the time it takes $i$ to have
  $|\hat{x}^i-x^i|$ evolve from $0$ to $\eps_i/(2\sqrt{\Bdout^i})$
  (note the conservativeness in this step as we disregard the first
  term on the righthand side
  of~\eqref{eq::TrigLaw_Distributed_UD}). The result of this analysis
  yields~\eqref{eq::state_holding_error}, but with
  %\begin{align*}
   % |\hat{x}^i-{x}^i (t)| \leq \frac{c^i}{\alpha}(\text{e}^{\alpha
    %  (t-t^i_{k})}-1), \quad t \ge t^i_{k} .
  %\end{align*}
  %where
  the value of $c^i$ given in the statement (i.e., the expression
  in~\eqref{eq::ci} with $\eta$ substituted by $\zeta$ as defined
  in~\eqref{eq::z2N_bound}).  Therefore, the time it takes
  $|\hat{x}^i-{x}^i|$ to reach $ \eps_i/(2\sqrt{ \Bdout^i})$ is lower
  bounded by $\tau^i>0$ as given
  in~\eqref{eq::dist_event_time_lw_bnd_undirected}.  This fact also
  implies $T = \infty$. Finally,~\eqref{eq:auxx} together
  with~\eqref{eq::a_to_c_Solution} and~\eqref{eq::z2N_bound}, imply
  that for $i\in{\VV}$,
  \begin{equation*}
    \!|x^i(t)-\avrg{\mathsf{r}}(t)|\!\leq\!
    \sqrt{\zeta^2 \!+\!
      |\vect{r}^\top(\vect{x}(0)-\Bvect{\mathsf{r}}(0))|^2},
  \end{equation*}
  for $ t \in \realnonnegative$.
  % Next, we examine the limiting behavior of the response of the
  % algorithm~\eqref{eq::EventTrig_Alg} to establish a tighter bound
  % on
  % the limiting value of $|x^i-\avrg{\mathsf{r}}|$ as given
  % in~\eqref{eq::Alg_D_Event_ultimate_bound_UD}.
  Moreover, since $T=\infty$, from~\eqref{eq::dV_intermidate}, we have
  \begin{align*}
    \dot{V}\leq &
    -\frac{1}{2}\beta\Ulambda_2(1-\theta)\vect{z}_{2:N}^\top\vect{z}_{2:N}+
    \frac{\beta\Ulambda_2}{2}\bar{r}(t),
  \end{align*}
  for $t\!\in\!\realnonnegative$, where $\bar{r}(t) =
  \frac{2}{\beta\Ulambda_2}\|\vect{z}_{2:N}\|\|\vect{q}_{2:N}(0)\|e^{-\alpha
    t}-\theta\vect{z}_{2:N}^\top\vect{z}_{2:N} +
  \frac{1}{2\Ulambda_2}\|\vect{\eps}\|^2 +
  \frac{2\gamma}{\beta\Ulambda_2}\|\vect{z}_{2:N}\|$.  Note that
  $\bar{r}(t)<0$ for
  \begin{align*}%\label{eq::r_sign_neg}
    \|\vect{z}_{2:N}\|\geq & \frac{\|\vect{q}_{2:N}(0)\|e^{-\alpha
        t}}{\beta\Ulambda_2\theta}+\frac{\gamma}{\beta\Ulambda_2\theta}
    \\
    & + \sqrt{\!\Big(\frac{\|\vect{q}_{2:N}(0)\|e^{-\alpha
          t}}{\beta\Ulambda_2\theta}\!+\!\frac{\gamma}{\beta\Ulambda_2\theta}\Big)^2\!
      +\!
      \frac{\|\vect{\eps}\|^2}{2\Ulambda_2\theta}}\!=\!\bar{\zeta}(t).
  \end{align*}
  Therefore, for $t\!\in\!\realnonnegative$, as long as $\|
  \vect{z}_{2:N} (t)\| \geq \bar{\zeta}(t)$,
  \begin{align*}
    \dot{V}\leq &
    -\frac{1}{2}\beta\Ulambda_2(1-\theta)\vect{z}_{2:N}^\top\vect{z}_{2:N} .
  \end{align*}
  As a result, 
  \begin{equation*}
    \limsup_{t\to\infty} \|\vect{z}_{2:N} (t) \| \leq
    \frac{\gamma}{\beta\Ulambda_2} +
    \sqrt{\Big(\frac{\gamma}{\beta\Ulambda_2}\Big)^2 
      + \frac{\|\vect{\eps}\|^2}{2\Ulambda_2}}.
  \end{equation*}
%    \begin{equation*}
 %   \limsup_{t\to\infty} \|\vect{z}_{2:N} (t) \| \leq
 %   \frac{\gamma}{\beta\Ulambda_2\theta} +
 %   \sqrt{\Big(\frac{\gamma}{\beta\Ulambda_2\theta}\Big)^2 
  %    + \frac{\|\vect{\eps}\|^2}{2\Ulambda_2\theta}}.
  %\end{equation*}
  %
Here, we used $\theta\to 1$.
  On the other hand, $\lim_{t\to\infty} z_1(t)=0$.  Combining these
  facts with~\eqref{eq:auxx}
  yields~\eqref{eq::Alg_D_Event_ultimate_bound_UD}, concluding the
  proof. \boxend
\end{pf}

We should point out that, the guaranteed lower
bound~\eqref{eq::dist_event_time_lw_bnd_undirected} on the
inter-event-times is more conservative than strictly necessary. 
This is because in our development, in order to decouple the analysis of the lower bound on the  inter-event times of each agent from its out-neighbors', we have neglected the effect of the term $\frac{1}{4 \dout^i(t) }
\sum\nolimits_{j=1}^N {\mathsf{a}}_{ij}(t)|\hat{x}^i(t) - \hat{x}^j(t)|^2$
in~\eqref{eq::TrigLaw_Distributed_UD}. The simulations of Section~\ref{sec::num} show the implementation
of~\eqref{eq::TrigLaw_Distributed_UD} resulting in inter-event times
longer than the ones of the event-triggered
law~\eqref{eq::TrigLaw_Distributed_Directed}.

%We should point out that, as  observed in the proof of
%Proposition~\ref{prop::Alg_D_Event_UD}, the guaranteed lower
%bound~\eqref{eq::dist_event_time_lw_bnd_undirected} on the
%inter-event-times is more conservative than strictly necessary because
%we have neglected the effect of the term $\frac{1}{4 \dout^i(t) }
%\sum_{j=1}^N {\mathsf{a}}_{ij}(t)|\hat{x}^i(t) - \hat{x}^j(t)|^2$
%in~\eqref{eq::TrigLaw_Distributed_UD} to decouple the analysis.
%The simulations of Section~\ref{sec::num} show the implementation
%of~\eqref{eq::TrigLaw_Distributed_UD} resulting in inter-event times
%longer than the ones of the event-triggered
%law~\eqref{eq::TrigLaw_Distributed_Directed}.

\begin{figure}[t!]
  \unitlength=0.5in \centering
 %   \captionsetup[subfloat]{captionskip=-1pt}
  %%%%%%%%%%%%%%%%%%%%%% 1 %%%%%%%%%%%%%%%%%%%%%%%%%%%%%%%%%%%
  \subfloat{
    \includegraphics[trim=0mm 0mm 0mm 0mm,clip,width=3.36in, height=1.45in]{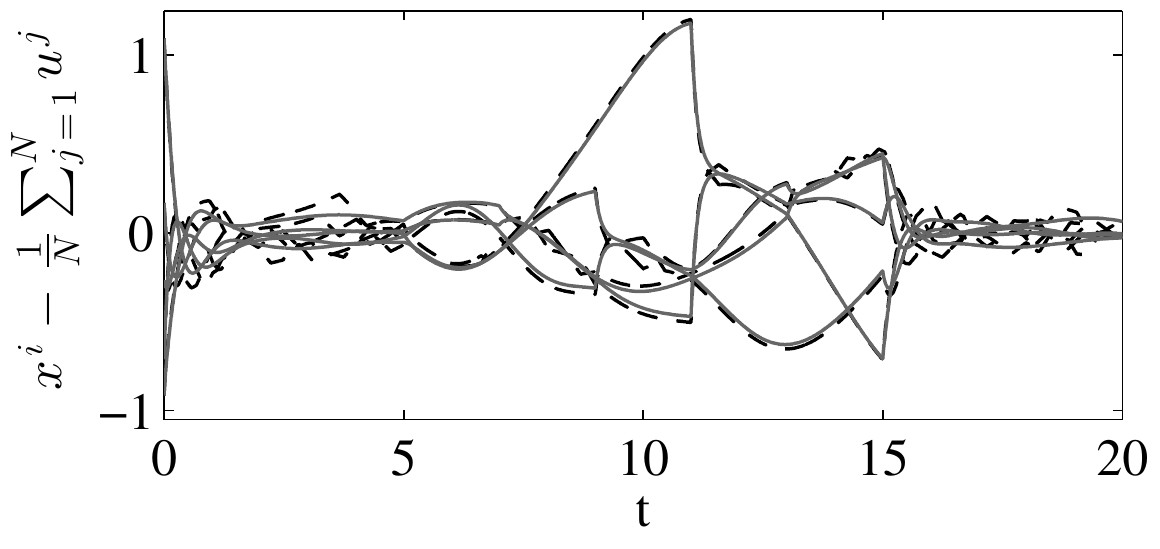}\quad\quad\quad
  }\\\vspace{-0.15in}
  %%%%%%%%%%%%%%%%%%%%%% 4 %%%%%%%%%%%%%%%%%%%%%%%%%%%%%%%%%%%
  \subfloat{\quad\,
    \includegraphics[width=3.2in, height=1.6in]{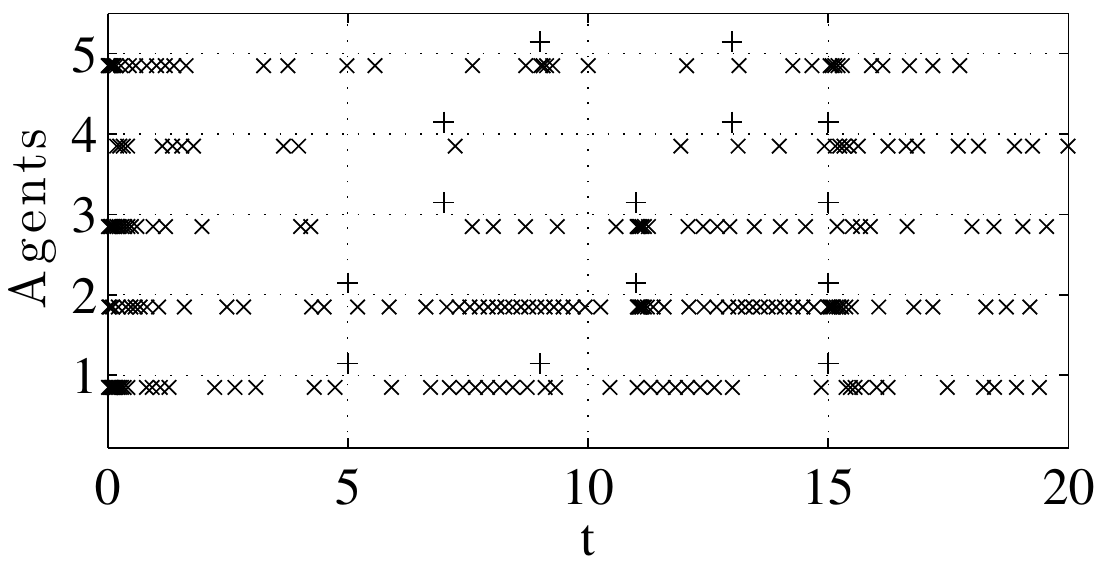} % }
  }\vspace{-0.15in}
  \caption{Executions of~\eqref{eq::EventTrig_Alg} with the
    event-triggered communication law described in
    Section~\ref{se:problem-statement} whose sampling rule
    is~\eqref{eq::TrigLaw_Distributed_Directed} and of~\eqref{eq::CT}
    with continuous-time communication.  The network is a
    weight-balanced time-varying digraph of $5$ agents with unit
    weights, where for $t\in[0,5)$ it is a fixed ring digraph, for
    $t\in[5,15)$ every $2$ seconds it is a single connected pair of
    nodes and then for $t\in[15,\infty)$ it is a ring digraph
    again. The inputs are $u^1(t)\!=\!0.5\sin(0.8t)$,
    $u^2(t)\!=\!0.5\sin(0.7t)\!+\!0.5\cos(0.6t)$,
    $u^3(t)=\sin(0.2t)\!+\!1$, $u^4(t)\!=\!\text{atan}(0.5t)$,
    $u^5(t)\!=\!0.1\cos(2t)$.  The top plot shows the tracking error
    and the bottom one shows the communication times of each agent.
    Black dashed (resp. gray solid) lines correspond to the
    event-triggered strategy with the sampling
    law~\eqref{eq::TrigLaw_Distributed_Directed} with $\eps_i\!=\!0.1$
    (resp. continuous-time communication~\eqref{eq::CT}).  In both
    cases, $\alpha\!=\!1$ and $\beta\!=\!4$. For each agent, $\times$
    indicates a sampling (and broadcast) time and $+$ indicates a time
    when the agent acquires a new in-neighbor.}\label{fig::Ex1sim}
\end{figure}

\section{Simulations}\label{sec::num}

\begin{figure}[t!]
  \unitlength=0.5in \centering
  % \captionsetup[subfloat]{captionskip=-1pt}
  %%%%%%%%%%%%%%%%%%%%%% 1 %%%%%%%%%%%%%%%%%%%%%%%%%%%%%%%%%%%
  \subfloat{
    \includegraphics[trim=0mm 0mm 0mm 0mm,clip,width=3.36in, height=1.4in]{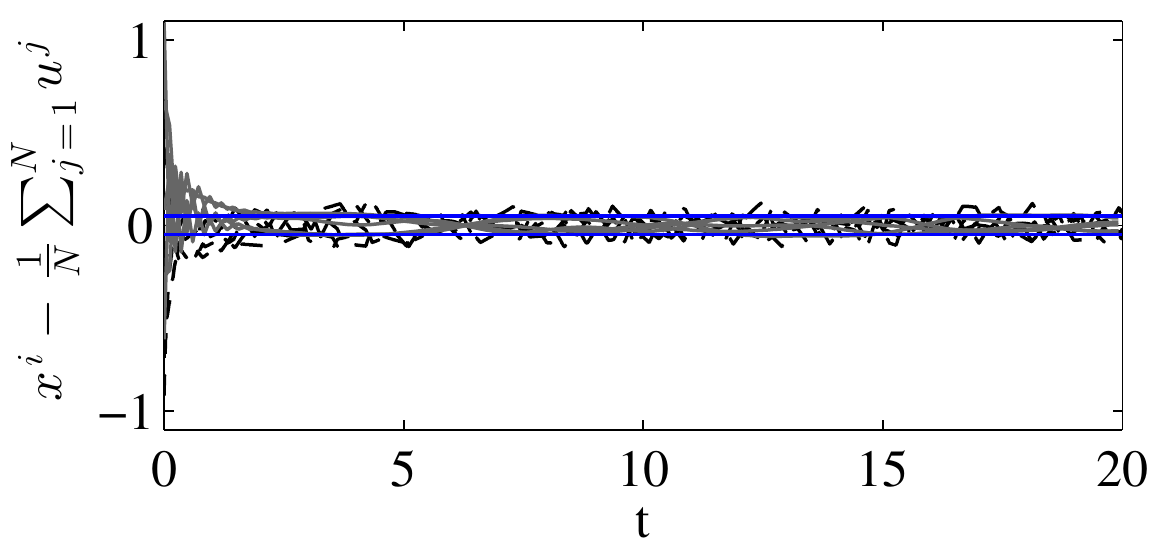}\quad\quad\quad
  }\\  \vspace{-0.15in}
  \subfloat{
    \includegraphics[trim=0mm 0mm 0mm 0mm,clip,width=3.36in, height=1.4in]{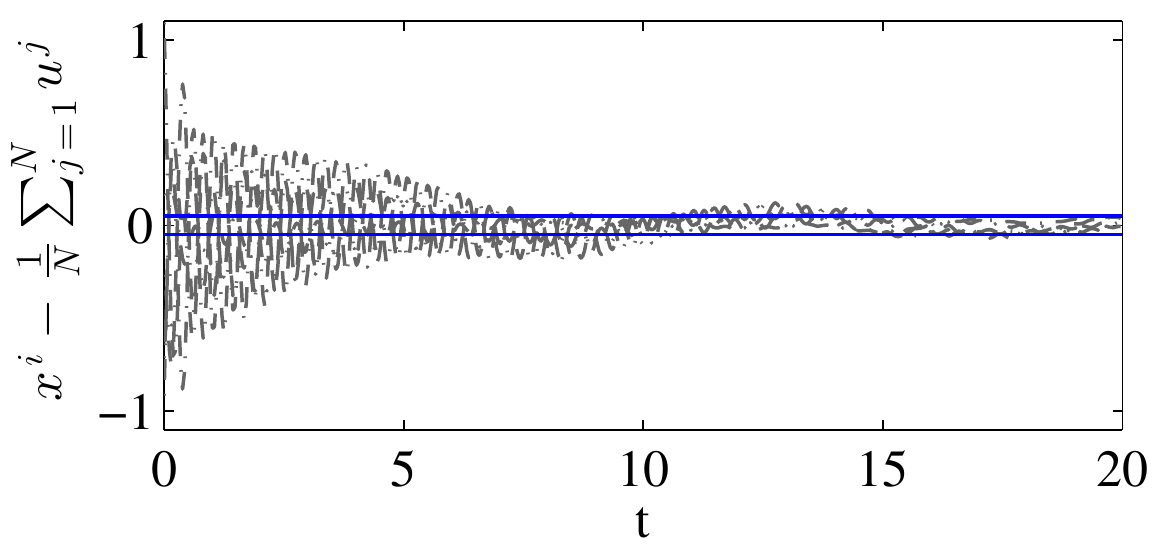}\quad\quad\quad
  }\\ \vspace{-0.15in}
  %%%%%%%%%%%%%%%%%%%%%% 4 %%%%%%%%%%%%%%%%%%%%%%%%%%%%%%%%%%%
  \subfloat{\quad\,
    \includegraphics[width=3.2in, height=1.4in]{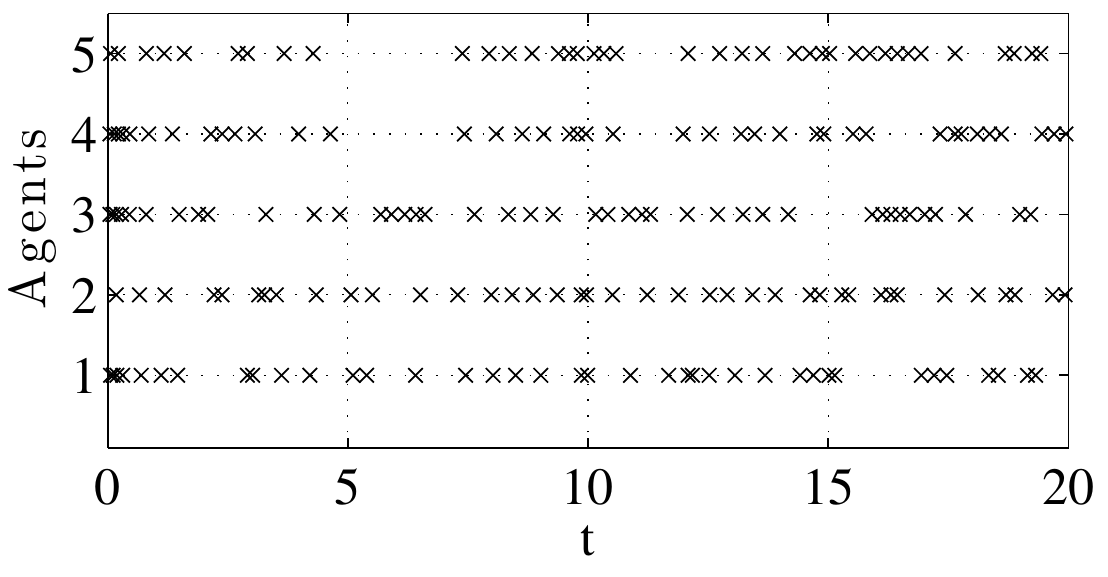} % }
  }\vspace{-0.15in}
  \caption{Comparison between the algorithm~\eqref{eq::EventTrig_Alg}
    employing the event-triggered communication law described in
    Section~\ref{se:problem-statement} with sampling
    rule~\eqref{eq::TrigLaw_Distributed_UD} and the Euler
    discretizations of the algorithm~\eqref{eq::CT} and the
    proportional-integral (PI) dynamic average consensus algorithm
    proposed in~\citep{RAF-PY-KML:06}. For the first two, we set
    $\alpha=1$ and $\beta=4$. For the latter, we set $\gamma=5$,
    $\lL_{\text{P}}=\lL$ and $\lL_{\text{I}}=4\lL$.  The network is a
    connected ring graph of $5$ agents with unit weights and the
    inputs are the same of Figure~\ref{fig::Ex1sim}.  In the top plot,
    black (resp. gray) lines correspond to the event-triggered
    law~\eqref{eq::TrigLaw_Distributed_UD} with
    $\eps_i/(2\sqrt{\dout^i})=0.1$ (resp. the Euler discretization of
    the algorithm~\eqref{eq::CT} with fixed stepsize
    $\delta=0.12$). The middle plot shows the response of the Euler
    discretization of the PI algorithm with fixed stepsize
    $\delta=0.039$. The horizontal lines in both the top and middle
    plots show the $\pm 0.05$ error bound for reference.  The bottom
    plot shows the communication times of each agent using the
    event-triggered strategy.}\label{fig::Ex2sim}
\end{figure}

\begin{figure}[t!]
  \unitlength=0.5in \centering
  % \captionsetup[subfloat]{captionskip=-1pt}
  %%%%%%%%%%%%%%%%%%%%%% 1 %%%%%%%%%%%%%%%%%%%%%%%%%%%%%%%%%%%
  \subfloat{
    \includegraphics[trim=1mm 0mm 1mm 0mm,clip,width=3.35in, height=1.4in]{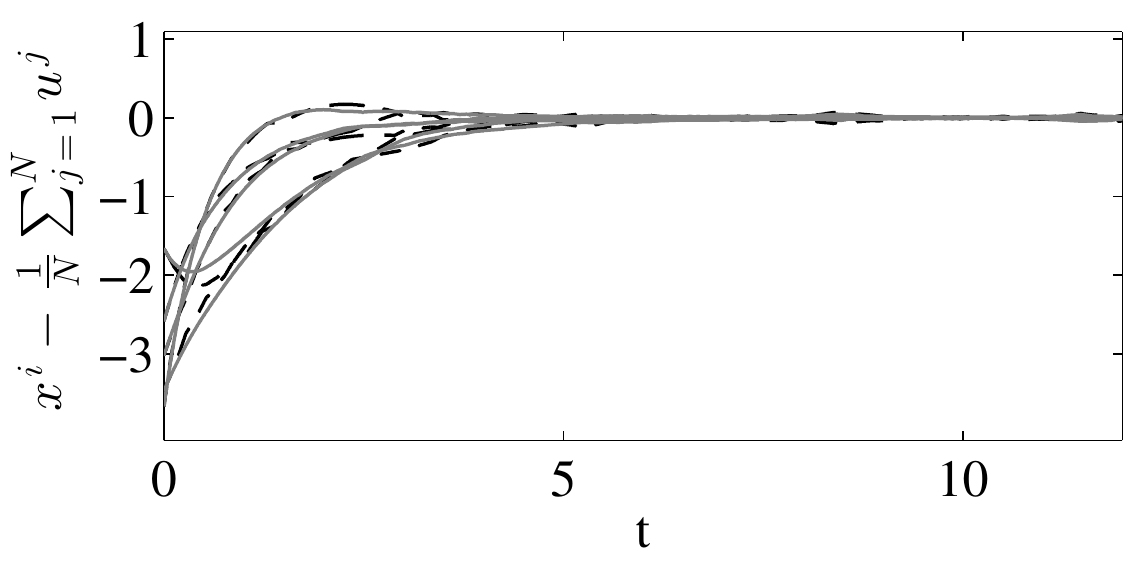}\quad\quad\quad
  }\\ \vspace{-0.15in}
  %%%%%%%%%%%%%%%%%%%%%% 4 %%%%%%%%%%%%%%%%%%%%%%%%%%%%%%%%%%%
  \subfloat{\quad\,
    \includegraphics[trim=1mm 0mm 1mm 0mm,clip,width=3.2in, height=1.6in]{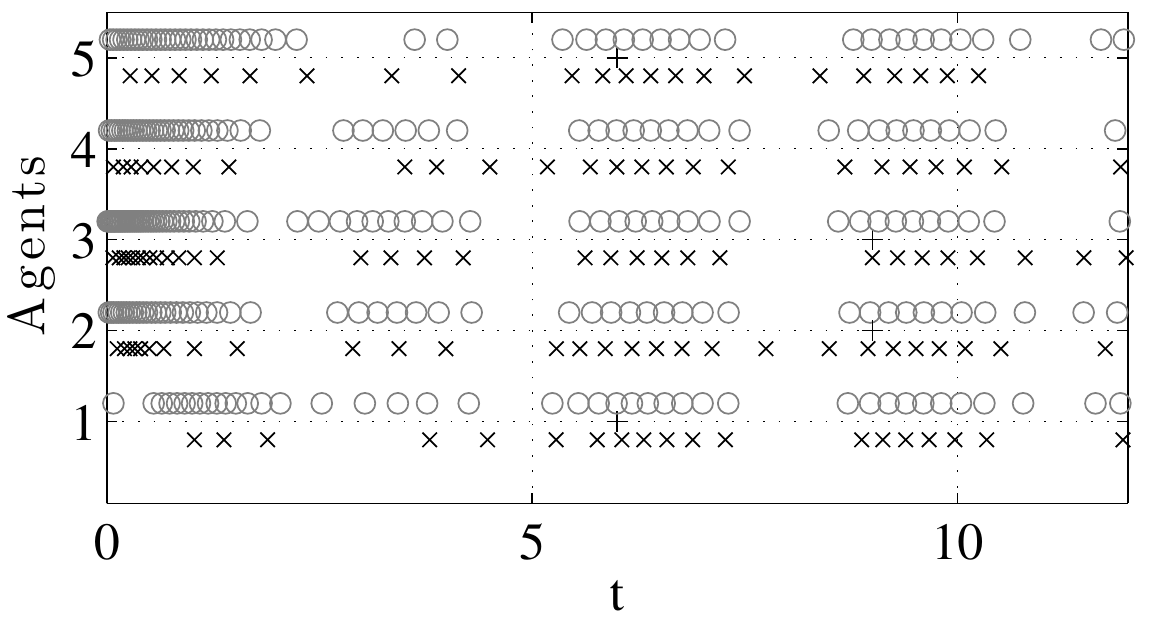} % }
  }\vspace{-0.15in}
  \caption{Executions of~\eqref{eq::EventTrig_Alg} with the sampling
    rule~\eqref{eq::TrigLaw_Distributed_Directed}
    and~\eqref{eq::TrigLaw_Distributed_UD} with the event-triggered
    communication law described in Section~\ref{se:problem-statement}.
    The network is a time-varying graph of $5$ agents corresponding to
    a connected ring graph with unit weights where one edge breaks
    every 3 seconds. The inputs are $u^1(t)\!=\!0.5\sin(t)+1/(t+2)+2$,
    $u^2(t)\!=\!0.5\sin(t)+1/(t+2)^2+4$,
    $u^3(t)\!=\!0.5\sin(t)+1/(t+2)^3+5$,
    $u^4(t)\!=\!0.5\sin(t)+\text{e}^{-t}+4$,
    $u^5(t)\!=\!0.5\sin(t)+\text{atan}(t)-1.5$. The top plot shows the
    tracking error with the gray solid (resp. black dashed) lines
    correspond to the law~\eqref{eq::TrigLaw_Distributed_Directed}
    with $\eps_i\!=\!0.1$ (resp. \eqref{eq::TrigLaw_Distributed_UD}
    with $\eps_i/(2 \sqrt{\dout^i})\!=\!0.1$). In both cases,
    $\alpha=\beta=1$. The bottom plot shows the communication times of
    each agent with o (resp. $\times$) markers corresponding to the
    law~\eqref{eq::TrigLaw_Distributed_Directed}
    (resp. \eqref{eq::TrigLaw_Distributed_UD}). For both cases, $+$
    shows the broadcasts associated to the acquisition of new
    in-neighbors. }\label{fig::Ex3sim}
\end{figure}

In this section, we illustrate the performance of the coordination
algorithm~\eqref{eq::EventTrig_Alg} under the event-triggered
communication laws~\eqref{eq::TrigLaw_Distributed_Directed}
and~\eqref{eq::TrigLaw_Distributed_UD} over a recurrently jointly
strongly connected digraph (cf. Fig.~\ref{fig::Ex1sim}),
a ring graph (cf. Fig.~\ref{fig::Ex2sim}) and a time-varying connected
graph (cf. Fig.~\ref{fig::Ex3sim}).  Figure~\ref{fig::Ex1sim} shows a
small degradation between the tracking performance of the
algorithm~\eqref{eq::EventTrig_Alg} with the event-triggered
communication law~\eqref{eq::TrigLaw_Distributed_Directed} and the
algorithm~\eqref{eq::CT} with continuous-time communication.  In the
event-triggered implementation, the number of times that agents
$\{1,2,3,4,5\}$ communicate in the time interval $[0,20]$ is
$(41,49,44,31,40)$, respectively. The large error observed in the time
interval $[5,15)$ is expected as in this time period every two seconds
only two agents are communicating with each other.  Naturally, these
two agents tend to converge to the average of their inputs and the
rest of the agents, being oblivious to the inputs of the other agents,
follow their own input.

Figure~\ref{fig::Ex2sim} compares the
algorithm~\eqref{eq::EventTrig_Alg} with event-triggered
communication~\eqref{eq::TrigLaw_Distributed_UD} and the Euler
discretizations of the algorithm~\eqref{eq::CT} and the
proportional-integral (PI) dynamic average consensus algorithm
proposed in~\citep{RAF-PY-KML:06}.  We set the parameters of the PI
algorithm so that its ultimate tracking error is similar to that
of~\eqref{eq::CT}.  For the discretizations, we use the largest
possible fixed stepsize $\delta=0.039$ for the PI algorithm (beyond
this value the algorithm diverges) and we use the stepsize
$\delta=0.12$ for the algorithm~\eqref{eq::CT}
(from~\citep{SSK-JC-SM:14-ijrnc}, convergence is guaranteed if $\delta
\in
(0,\min\{\alpha^{-1},\beta^{-1}(\text{d}_{\max}^{\text{out}})^{-1}\})$,
which for this example results in $\delta\in(0,0.125)$).  The number
of times that agents $\{1,2,3,4,5\}$ communicate in the time interval
$[0,20]$ is $(39, 40,42, 40,39)$, respectively, when implementing
event-triggered communication~\eqref{eq::TrigLaw_Distributed_UD}. This
is significantly less than the communication used by each agent in the
Euler discretizations of~\eqref{eq::CT} ($20/0.12 \simeq 166$ rounds)
and the PI algorithm ($20/0.039 \simeq 512$ rounds).
   
Figure~\ref{fig::Ex3sim} shows the execution
of~\eqref{eq::EventTrig_Alg} with the event-triggered communication
laws~\eqref{eq::TrigLaw_Distributed_Directed}
and~\eqref{eq::TrigLaw_Distributed_UD} over a time-varying connected
graph. For each agent $i \in \{1,2,3,4,5\}$, we choose $\eps_i$ for
each law so that the summand in the right-hand side of the trigger
($\eps_i$ for~\eqref{eq::TrigLaw_Distributed_Directed}, $\eps_i/(2
\sqrt{\Bdout^i})$ for~\eqref{eq::TrigLaw_Distributed_UD}) amounts to
the same quantity.  The plots show similar tracking performance for
both algorithms, with the law~\eqref{eq::TrigLaw_Distributed_UD}
inducing less than half communication
than~\eqref{eq::TrigLaw_Distributed_Directed}. In fact, the number of
times that agents $\{1,2,3,4,5\}$ communicate in the time interval
$[0,12]$ is $(43,56,72,57,55)$
under~\eqref{eq::TrigLaw_Distributed_Directed} and $(21,30,33,26,23)$
under~\eqref{eq::TrigLaw_Distributed_UD}.

\section{Conclusions}\label{sec::conclu}
We have studied the multi-agent dynamic average consensus problem over
networks where inter-agent communication takes place at discrete time
instants in an opportunistic fashion.  Our starting point has been our
previously developed continuous-time dynamic average consensus
algorithm which is known to converge exponentially to a small
neighborhood of the network's inputs average. We have proposed two
different distributed event-triggered laws that agents can employ to
trigger communication with neighbors, depending on whether the
interaction topology is described by a weight-balanced and recurrently
jointly strongly connected digraph or a time-varying connected
undirected graph. In both cases, we have established the correctness
of the algorithm and showed that a positive lower bound on the
inter-event times of each agent exists, ruling out the presence of
Zeno behavior.  Future work will be devoted to further relaxing the
connectivity requirements on the interaction topology (tying them in
to the evolution of the dynamic inputs available to the agents), the
improvement of the practical convergence guarantees using time-varying
thresholds in the trigger design, the use of agent abstractions in the
development of self-triggered communication laws, and the synthesis of
other distributed triggers that individual agents can evaluate
autonomously and lead to a more efficient use of the limited network
resources.

%% FOR FINAL VERSION
% \section*{Acknowledgements}
% The work of the first author is supported by a UC President's
% Postdoctoral Fellowship and the Cymer Center for Control Systems and
% Dynamics. This research was partially supported by NSF award
% ECCS-1307176.

{\footnotesize%

}


\begin{thebibliography}{22}
\providecommand{\natexlab}[1]{#1}
\providecommand{\url}[1]{\texttt{#1}}
\expandafter\ifx\csname urlstyle\endcsname\relax
  \providecommand{\doi}[1]{doi: #1}\else
  \providecommand{\doi}{doi: \begingroup \urlstyle{rm}\Url}\fi

\bibitem[Bai et~al.(2010)Bai, Freeman, and Lynch]{HB-RAF-KML:10}
H.~Bai, R.~A. Freeman, and K.~M. Lynch.
\newblock Robust dynamic average consensus of time-varying inputs.
\newblock In \emph{{IEEE} Int. Conf. on Decision and Control}, pages
  3104--3109, Atlanta, GA, USA, December 2010.

\bibitem[Bullo et~al.(2009)Bullo, Cort{\'e}s, and Mart{\'\i}nez]{FB-JC-SM:09}
F.~Bullo, J.~Cort{\'e}s, and S.~Mart{\'\i}nez.
\newblock \emph{Distributed Control of Robotic Networks}.
\newblock Applied Mathematics Series. Princeton University Press, 2009.
\newblock ISBN 978-0-691-14195-4.
\newblock Available at http://www.coordinationbook.info.

\bibitem[Carron et~al.(2013)Carron, Todescato, Carli, and
  Schenato]{AC-TO-MT-RC-LS:13}
A.~Carron, M.~Todescato, R.~Carli, and L.~Schenato.
\newblock Adaptive consensus-based algorithms for fast estimation from relative
  measurements.
\newblock In \emph{{IFAC} Workshop on Distributed Estimation and Control in
  Networked Systems}, Koblenz, Germany, March 2013.

\bibitem[Dimarogonas et~al.(2012)Dimarogonas, Frazzoli, and
  Johansson]{DVD-EF-KHJ:12}
D.~V. Dimarogonas, E.~Frazzoli, and K.~H. Johansson.
\newblock Distributed event-triggered control for multi-agent systems.
\newblock \emph{IEEE Transactions on Automatic Control}, 57\penalty0
  (5):\penalty0 1291--1297, 2012.

\bibitem[Donkers and Heemels(2012)]{MCFD-WPMHH:12}
M.~C.~F. Donkers and W.~P. M.~H. Heemels.
\newblock Output-based event-triggered control with guaranteed
  $\mathcal{L}_\infty$-gain and improved and decentralized event-triggering.
\newblock \emph{IEEE Transactions on Automatic Control}, 57\penalty0
  (6):\penalty0 1362--1376, 2012.

\bibitem[Fan et~al.(2013)Fan, Feng, Wang, and Song]{YF-GF-YW-CS:13}
Y.~Fan, G.~Feng, Y~Wang, and C.~Song.
\newblock Distributed event-triggered control of multi-agent systems with
  combinational measurements.
\newblock \emph{Automatica}, 49\penalty0 (2):\penalty0 671--675, 2013.

\bibitem[Freeman et~al.(2006)Freeman, Yang, and Lynch]{RAF-PY-KML:06}
R.~A. Freeman, P.~Yang, and K.~M. Lynch.
\newblock Stability and convergence properties of dynamic average consensus
  estimators.
\newblock In \emph{{IEEE} Int. Conf. on Decision and Control}, pages 338--343,
  2006.

\bibitem[Garcia et~al.(2013)Garcia, Cao, Yuc, Antsaklis, and
  Casbeer]{EG-YC-HY-PA-DC:13}
E.~Garcia, Y.~Cao, H.~Yuc, P.~Antsaklis, and D.~Casbeer.
\newblock Decentralised event-triggered cooperative control with limited
  communication.
\newblock \emph{International Journal of Control}, 86\penalty0 (9):\penalty0
  1479--1488, 2013.

\bibitem[Heemels et~al.(2012)Heemels, Johansson, and Tabuada]{WPMHH-KHJ-PT:12}
W.~P. M.~H. Heemels, K.H Johansson, and P.~Tabuada.
\newblock An introduction to event-triggered and self-triggered control.
\newblock In \emph{{IEEE} Int. Conf. on Decision and Control}, pages
  3270--3285, Maui, HI, 2012.

\bibitem[Khalil(2002)]{HKK:02}
H.~Khalil.
\newblock \emph{Nonlinear Systems}.
\newblock Prentice Hall, 2002.

\bibitem[Kia et~al.(2014{\natexlab{a}})Kia, Cort{\'e}s, and
  Mart{\'\i}nez]{SSK-JC-SM:14-cdc}
S.~S. Kia, J.~Cort{\'e}s, and S.~Mart{\'\i}nez.
\newblock Dynamic average consensus with distributed event-triggered
  communication.
\newblock In \emph{{IEEE} Int. Conf. on Decision and Control}, Los Angeles, CA,
  2014{\natexlab{a}}.
\newblock Submitted.

\bibitem[Kia et~al.(2014{\natexlab{b}})Kia, Cort\'es, and
  Mart{\'\i}nez]{SSK-JC-SM:14-ijrnc}
S.~S. Kia, J.~Cort\'es, and S.~Mart{\'\i}nez.
\newblock Dynamic average consensus under limited control authority and privacy
  requirements.
\newblock \emph{International Journal on Robust and Nonlinear Control},
  2014{\natexlab{b}}.
\newblock To appear.

\bibitem[Mazo and Tabuada(2011)]{MMJ-PT:11}
M.~Mazo and P.~Tabuada.
\newblock Decentralized event-triggered control over wireless sensor/actuator
  networks.
\newblock \emph{IEEE Transactions on Automatic Control}, 56\penalty0
  (10):\penalty0 2456--2461, 2011.

\bibitem[Nowzari and Cort{\'e}s(2014)]{CN-JC:14-acc}
C.~Nowzari and J.~Cort{\'e}s.
\newblock Zeno-free, distributed event-triggered communication and control for
  multi-agent average consensus.
\newblock In \emph{{A}merican {C}ontrol {C}onference}, Portland, OR, 2014.
\newblock To appear.

\bibitem[Olfati-Saber(2007)]{ROS:07}
R.~Olfati-Saber.
\newblock Distributed {K}alman filtering for sensor networks.
\newblock In \emph{{IEEE} Int. Conf. on Decision and Control}, pages
  5492--5498, New Orleans, LA, December 2007.

\bibitem[Olfati-Saber and Shamma(2005)]{ROS-JSS:05}
R.~Olfati-Saber and J.~S. Shamma.
\newblock Consensus filters for sensor networks and distributed sensor fusion.
\newblock In \emph{{IEEE} Int. Conf. on Decision and Control and European
  Control Conference}, pages 6698--6703, Seville, Spain, December 2005.

\bibitem[Seyboth et~al.(2013)Seyboth, Dimarogonas, and
  Johansson]{GSS-DVD-KHJ:13}
G.~S. Seyboth, D.~V. Dimarogonas, and K.~H. Johansson.
\newblock Event-based broadcasting for multi-agent average consensus.
\newblock \emph{Automatica}, 49\penalty0 (1):\penalty0 245--252, 2013.

\bibitem[Spanos et~al.(2005)Spanos, Olfati-Saber, and Murray]{DPS-ROS-RMM:05b}
D.~P. Spanos, R.~Olfati-Saber, and R.~M. Murray.
\newblock Dynamic consensus on mobile networks.
\newblock In \emph{{IFAC} {W}orld {C}ongress}, Prague, Czech Republic, July
  2005.

\bibitem[Wang and Lemmon(2011)]{XW-MDL:11}
X.~Wang and M.~D. Lemmon.
\newblock Event-triggering in distributed networked control systems.
\newblock \emph{IEEE Transactions on Automatic Control}, 56\penalty0
  (3):\penalty0 586--601, 2011.

\bibitem[Yang et~al.(2007)Yang, Freeman, and Lynch]{PY-RAF-KML:07}
P.~Yang, R.~A. Freeman, and K.~M. Lynch.
\newblock Distributed cooperative active sensing using consensus filters.
\newblock In \emph{{IEEE} Int. Conf. on Robotics and Automation}, pages
  405--410, Roma, Italy, April 2007.

\bibitem[Yang et~al.(2008)Yang, Freeman, and Lynch]{PY-RAF-KML:08}
P.~Yang, R.~A. Freeman, and K.~M. Lynch.
\newblock Multi-agent coordination by decentralized estimation and control.
\newblock \emph{IEEE Transactions on Automatic Control}, 53\penalty0
  (11):\penalty0 2480?--2496, 2008.

\bibitem[Zhu and Mart{\'\i}nez(2010)]{MZ-SM:08a}
M.~Zhu and S.~Mart{\'\i}nez.
\newblock Discrete-time dynamic average consensus.
\newblock \emph{Automatica}, 46\penalty0 (2):\penalty0 322--329, 2010.

\end{thebibliography}
\end{document}